\documentclass[preprint,11pt]{elsarticle}
\usepackage{amsmath,amssymb,amsthm,mathrsfs,geometry,graphicx,setspace,soul}
\usepackage{multirow}
\usepackage[table]{xcolor}
\usepackage{hyperref}

\usepackage{subcaption}

\definecolor{aura}{RGB}{102,181,91}
\definecolor{darkRed}{RGB}{140,0,0}

\usepackage{pifont}
\newcommand{\cmark}{{\color{aura} \ding{51}}}%
\newcommand{\xmark}{{\color{darkRed} \ding{55}}}%

\usepackage{scalerel}
\DeclareMathOperator*{\A}{\scalerel*{A}{\sum}}

\usepackage{tikz}
\usepackage{collcell}

\newcommand*{\MinNumber}{0.2}%
\newcommand*{\MaxNumber}{2.0}%

\definecolor{minColour}{RGB}{255,255,0}	
\definecolor{maxColour}{RGB}{102,181,91}	

\newcommand{\ApplyGradient}[1]{%
	\ifdim #1 pt < \MinNumber pt
		\colorbox{minColour}{#1}
	\else
       		\ifdim #1 pt < \MaxNumber pt
            		\pgfmathsetmacro{\PercentColor}{max(min(100.0*(#1 - \MinNumber)/(\MaxNumber-\MinNumber),100.0),0.00)} %
           		 \hspace{-0.33em}\colorbox{maxColour!\PercentColor!minColour}{#1}
       		\else
           		\colorbox{maxColour}{#1}
        		\fi
	\fi
}

\newcolumntype{R}{>{\collectcell\ApplyGradient}c<{\endcollectcell}}
\usepackage{amssymb}

\biboptions{sort&compress,square}

\def\presup#1{\kern+.1em\relax{}^{#1}\kern-.1em\relax}      

\journal{arXiv}

\begin{document}

\begin{frontmatter}



\title{Ghost stabilisation of the Material Point Method for stable quasi-static and dynamic analysis of large deformation problems}


\author[du]{William M. Coombs}

\address[du]{Department of Engineering, Durham University, UK.\\
E: w.m.coombs@durham.ac.uk}

\begin{abstract}

  The unstable nature of the material point method is widely documented and is a barrier to the method being used for routine engineering analyses of large deformation problems.  The vast majority of papers concerning this issue are focused on the instabilities that manifest when a material point crosses between background grid cells.  However, there are other issues related to the stability of material point methods.  This paper focuses on the issue of the conditioning of the global system of equations caused by the arbitrary nature of the position of the physical domain relative to the background computational grid.  The issue is remedied here via the use of a Ghost stabilisation technique that penalises variations in the gradient of the solution field near the boundaries of the physical domain.  This technique transforms the stability of the material point method, providing a robust computational framework for large deformation dynamic and \emph{quasi}-static analysis.        

\end{abstract}


\end{frontmatter}


\section{Introduction}

The material point method's \cite{Sulsky1994} greatest advantage, the decoupling of the deformation of the physical material from the computational grid, is also its greatest challenge in terms of robust and efficient engineering computations.  By allowing the physical bodies to deform through a background grid of finite elements, the Material Point Method (MPM) can be used to model large deformation solid mechanics problems involving history dependent materials without re-meshing or re-mapping of material parameters (see \cite{Vaucorbeil2020,Solowski2021} for recent reviews articles).  However, this means the overlap between the physical body, represented by material points, and the background grid, which is used to solve the governing equations, varies over the analysis and this causes numerical stability problems \cite{Ma2010}. Specifically, background elements that are partially filled with material points cause problems in terms of the conditioning of the global system of equations.  Similar problems have been observed in the cut finite element method \cite{Sticko2020,Hansbo2017,Burman2018}, where arbitrarily small intersections can occur between an element, $K$, and the physical problem domain, $\Omega$.   
These intersections can be much smaller than the size of the elements, $h$.  This impacts on the conditioning of the linear system of equations being solved as the smallest eigenvalue of the linear system is related to the smallest intersection between the elements and the background grid, which can be arbitrarily small.  This means that the condition number (the ratio of the largest to smallest eigenvalue) of the  linear system of equations is not bounded.  These issues have a detrimental impact on both dynamic and \emph{quasi}-static analysis via the conditioning of the mass and/or stiffness matrix, with the issues often being more severe for the stiffness matrix \citep{Sticko2020}.  In terms of material point methods, the issue is compounded when using domain-based MPMs, such as the Generalised Interpolation Material Point Method (GIMPM) \cite{Bardenhagen2004}, where arbitrarily small overlaps can occur between the material point domains, $\Omega_p$, and the elements in the background mesh.
A crude way avoid this issue is to exclude \emph{small} overlaps between the material point domains and the background mesh (or rather \emph{small} shape function values).  However, the definition of \emph{small} is arbitrary and problem dependent and does not enforce a bound on the conditioning of the linear system of equations being solved.  

The conditioning issues described above means that the vast majority of explicit dynamic material point methods adopt a \emph{lumped} mass matrix, making the inversion of the mass matrix trivial.  However, adopting a lumped mass matrix has consequences for energy conservation \cite{Burgess1992,Bardenhagen2002,Love2006,Nairn2021,Pretti2022}, resulting in artificial numerical dissipation.  This issue was recognised by \citet{Nairn2021} who proposed an approximation to the consistent mass matrix inverse based on an eXtended Particle-In-Cell (XPIC) concept \cite{Hammerquist2017} which allowed the mass matrix inverse to be written as an infinite series expansion.  Small overlaps between material points and the background grid causes other issues.  For example, \citet{Ma2010} highlighted impact of small background grid nodal masses on the stability of material point methods when determining the grid node accelerations.  They introduced an approach to redistribute the nodal forces associated with the small nodal masses to eliminate the spurious accelerations.  However, the method does not deal with the small values in the mass matrix, choosing instead to deal with the manifestation of these small values in the acceleration calculation. Very few solutions to this issue have been proposed for implicit \emph{quasi}-static analyses, where the conditioning of the stiffness matrix is critical for stable convergence.  \citet{Wang2016} are proponents of a \emph{soft stiffness} to stabilisation, where it is assumed that the background finite element mesh has a \emph{small} elastic modulus, which is integrated and added to the global stiffness matrix.  However, the approach is not targeted at the regions that need stabilisation and degrades the equilibrium convergence rate of implicit solution methods.  

This paper resolves the issue of ill-conditioning of the global consistent mass and stiffness matrices for dynamic and \emph{quasi}-static material point methods using, and adapting, the \emph{Ghost stabilisation} technique of \citet{Burman2010}.  The method is general in that it can be applied to any material point variant and does not require the physical boundary of the domain to be tracked explicitly.  The Ghost stabilisation opens the door for the consistent mass matrix to be used for routine dynamic material point analyses, eliminating the additional energy dissipation caused by a lumped mass matrix and provides flexibility in the stress and velocity updating procedures.  It offers the first credible stabilisation technique for implicit \emph{quasi}-static material point methods, removing much of the uncertainty regarding if a given analysis will converge or not.  Overall, this paper provides a key ingredient for material point methods to become a practical choice for large deformation engineering analyses.


\section{Material Point Method}

\begin{figure}[!h]
\centering
    \begin{subfigure}[t]{0.33\textwidth}
        \centering
        \includegraphics[width=0.99\textwidth]{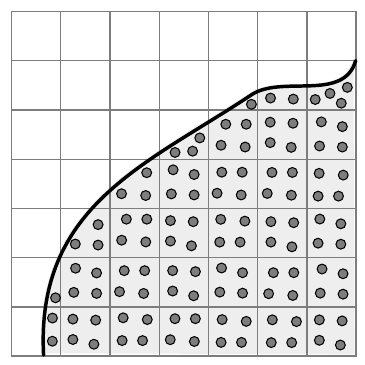}
        \caption{\footnotesize initial position}
    \end{subfigure}%
    \begin{subfigure}[t]{0.33\textwidth}
        \centering
        \includegraphics[width=0.99\textwidth]{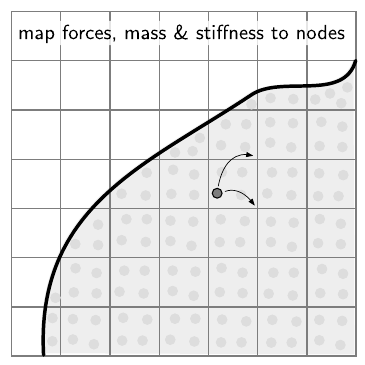}
        \caption{\footnotesize point-to-grid map}
    \end{subfigure}%
    \begin{subfigure}[t]{0.33\textwidth}
        \centering
        \includegraphics[width=0.99\textwidth]{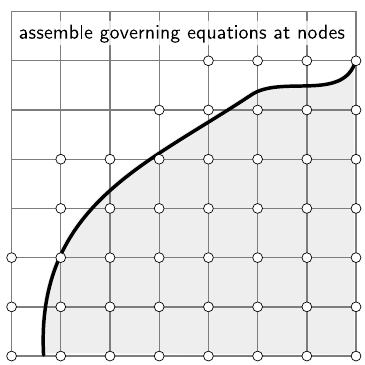}
        \caption{\footnotesize assembly}
    \end{subfigure} \newline
    \begin{subfigure}[t]{0.33\textwidth}
        \centering
        \includegraphics[width=0.99\textwidth]{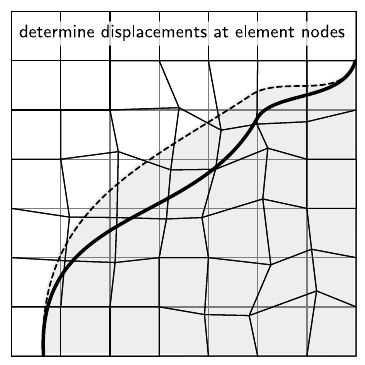}
        \caption{\footnotesize grid solution}
    \end{subfigure}%
    \begin{subfigure}[t]{0.33\textwidth}
        \centering
        \includegraphics[width=0.99\textwidth]{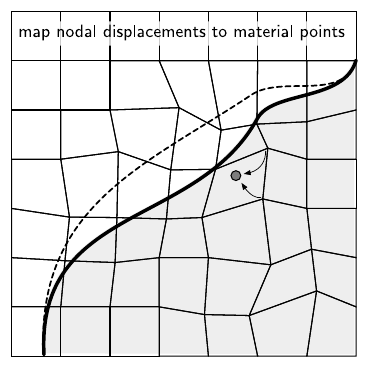}
        \caption{\footnotesize grid-to-point map}
    \end{subfigure}%
    \begin{subfigure}[t]{0.33\textwidth}
        \centering
        \includegraphics[width=0.99\textwidth]{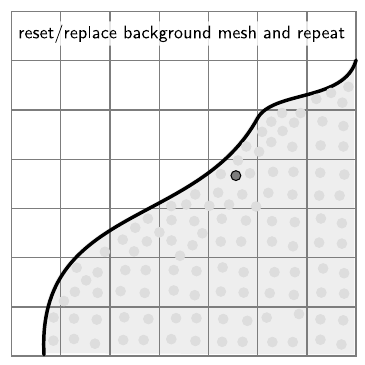}
        \caption{\footnotesize update points \& reset grid}
    \end{subfigure}
\caption{Material point method steps (adapted from \cite{Coombs2020}).}
\label{Fig:MPMsteps}
\end{figure}

The material point method uses material points to represent the physical body/bodies being analysed but the governing equations are solved on a background mesh/grid of finite elements.  This means that quantities must be mapped from points-to-grid and from grid-to-points during different stages of an analysis, and the key steps in a material point method algorithm are shown in Figure~\ref{Fig:MPMsteps}.  A given problem is discretised into a number of time (or load) steps and for each of the steps the following procedure is applied:
 \begin{enumerate}[(a)]
  \item the physical body is discretised by a number of material points that carry information on the volume, mass, deformation, stress, etc. of the material that they represent and the position of these points on the background grid determined (i.e. identifying the element(s)/nodes(s) the points interact with);
  \item material point information on velocity/momentum, stress, mass and stiffness is mapped to the background grid using the basis functions (or spatial derivatives) of the MPM variant under consideration to form nodal velocity, force, mass and stiffness quantities, as required by the adopted time stepping algorithm; 
  \item governing equations are assembled on the background grid nodes;
  \item governing equations are solved for the unknown nodal displacements/velocities (depending on the nature of the governing equations and the adopted time stepping algorithm);
 \item grid information is mapped to the material points to update deformation, velocity, stress, volume, etc. and to determine the material point displacement over the time step; 
 \item material points are updated\footnote{In the MPM literature this stage is often referred to as the \emph{convection} step, or that the material points are \emph{convected} to their new positions.  \emph{Convected} is an erroneous choice of word to describe this stage of the analysis as it has nothing to do with convection (especially in solid materials), it is simply saying that the material point positions are updated based on the their displacements within the current time/load step.} to their new positions and the background grid reset/replaced.  
\end{enumerate}
The above general steps are applicable to dynamic and \emph{quasi}-static material point methods discretised in time via implicit and explicit approaches.  However, specific details within each of the steps will change depending on algorithm/implementation choices, for example the material point stress updating algorithmic location will change if an Update Stress First (USF) explicit approach is adopted, as the above position corresponds to an Update Stress Last (USL) method.  Key to this procedure is the definition of the governing equations that will be solved during the analysis.  This paper is focused on dynamic and \emph{quasi}-static analysis of elastic and elasto-plastic solids undergoing large deformations and the governing equations for this continuum behaviour will be covered in the next section.

\subsection{Continuum formulation}

The material point approach adopted in this paper is based on the open-source AMPLE (A Material Point Learning Environment) code \cite{Coombs2020}, which is a \emph{quasi}-static implicit finite deformation elasto-plastic material point method based on an updated Lagrangian formulation\footnote{AMPLE adopts the same formulation as implemented in the generalised interpolation approach of \citet{Charlton2017}.}, but extended to include inertia effects for the purpose of this paper.  The method is defined by the following weak statement of equilibrium 
\begin{equation}\label{eqn:weak}
  \int_{\varphi_t(\Omega)}\Bigl( \sigma_{ij}(\nabla_x \eta)_{ij} -(b_i-\rho\dot{v}_i)\eta_i \Bigr) \text{d}V - \int_{\varphi_t(\partial\Omega)}\bigl(t_i\eta_i \bigr) \text{d}s = 0.
\end{equation}
where $\dot{v}_i$ and $\varphi_t$ are the acceleration and motion of the material body with domain, $\Omega$, which is subjected to tractions, $t_i$, on the boundary of the domain (with surface, $s$), $\partial\Omega$, and body forces, $b_i$, acting over the volume, $V$ of the domain, which has a density, $\rho$, and leads to a Cauchy stress field, $\sigma_{ij}$, through the body. The weak form is derived in the current frame assuming a field of admissible virtual displacements, $\eta_i$.  In the case of \emph{quasi}-static analysis the acceleration of the body is assumed to be zero, which reduces (\ref{eqn:weak}) to that of the AMPLE framework \cite{Coombs2020}.  The surface traction term in (\ref{eqn:weak}) is neglected as the focus of this paper is on the stability of the material point method rather than the imposition of surface tractions within the method, which is a research question in its own right (for example see \cite{Bing2019,Remmerswaal2017}).  We also restrict the Dirichlet boundary conditions to coincide with the background grid as the imposition of general Dirichlet constraints is another active area of research, for example \cite{Chandra2021,Bing2019,Cortis2018}.

The large deformation elasto-plastic continuum mechanics formulation used in this paper is identical to that of \citet{Charlton2017} and \citet{Coombs2020} and readers are referred to those papers for details.  In brief, the formulation adopts a multiplicative decomposition of the deformation gradient into elastic and plastic components, with a linear relationship between logarithmic elastic strain and Kirchhoff stress.  This is combined with an exponential map of the plastic flow rule which allows conventional small strain plasticity algorithms to be used without modification \cite{Simo1992a}.  This powerful combination is widely used in large deformation finite element methods \cite{SouzaNeto} and several material point method implementations \cite{Charlton2017,Coombs2018,Coombs2020,Cortis2018,Wang2019,Coombs2020a,Wang2021}.  These numerical implementations require the weak statement of equilibrium (\ref{eqn:weak}) to be discretised in space and time.

\subsection{Discrete formulation}

Neglecting the traction term, the Galerkin form of the weak statement of equilibrium over each background grid element, $K$, can be obtained from (\ref{eqn:weak}) as
\begin{equation}\label{eqn:MPMweak}
  \int_{\varphi_t(K)}[\nabla_x S_{vp}]^{T}\{\sigma_p\} \text{d}V - \int_{\varphi_t(K)}[S_{vp}]^{T}\{b\} \text{d}V + \int_{\varphi_t(K)}\rho[S_{vp}]^{T}\{\dot{v}\} \text{d}V= \{0\},
\end{equation}
where  $[S_{vp}]$ contains the basis/shape functions that are used to transfer information between the material points and the background grid, $[\nabla_x S_{vp}]$ is he strain-displacement matrix containing derivatives of the basis functions with respect to the updated coordinates\footnote{$[\nabla_x S_{vp}]$ is essentially the same as the conventional strain-displacement $[B]$ matrix found in finite element literature and $[S_{vp}]$ is the equivalent of the shape function matrix, often denoted $[N]$.}.  Consistent with the majority of the published MPM literature, the subscripts $(\cdot)_v$ and $(\cdot)_p$ denote quantities associated with \emph{vertices} (nodes) of the background grid and material \emph{points}, respectively.    

In the material point method the physical domain is discretised by a material points, each representing a finite volume of material, $V_p$.  This allows the Galerkin form of the weak statement of equilibrium (\ref{eqn:MPMweak}) to be approximated as 
\begin{equation}\label{eqn:MPMweakDiscrete}
  \A_{ \forall p} \Bigl( [\nabla_x S_{vp}]^{T}\{\sigma_p\} V_p - [S_{vp}]^{T}\{b\} V_p + [S_{vp}]^{T}\{\dot{v}\} m_p \Bigr)= \{0\},
\end{equation}
where $\A$ is the standard assembly operator acting over all of the material points in the problem and $m_p=\rho V_p$ is the mass associated with a material point.  The first term in (\ref{eqn:MPMweakDiscrete}) represents the internal nodal forces generated by the stress within the material, the second term is the external actions caused by the body forces and the third term accounts for the inertia of the physical body.

\subsection{Explicit dynamics}\label{sec:ExpDyn}

Explicit time stepping approaches assumed that the next state can be determined from information known at the current state, essentially forward predicting the new accelerations, $\{\dot{v}_v\}$, via the difference between the internal and external forces.    Assuming that the material point acceleration can be approximated as
\begin{equation}
	\{\dot{v}\} = [S_{vp}]\{\dot{v}_v\},
\end{equation}
where $\{\dot{v}_v\}$ are the accelerations of the background grid nodes, (\ref{eqn:MPMweakDiscrete}) becomes
\begin{equation}
  \A_{ \forall p} \Bigl( [\nabla_x S_{vp}]^{T}\{\sigma_p\} V_p - [S_{vp}]^{T}\{b\} V_p + [S_{vp}]^{T}[S_{vp}]\{\dot{v}_v\} m_p \Bigr)= \{0\},
\end{equation}
which can be rearranged to obtain the grid node accelerations 
\begin{equation}\label{eqn:nodeAcc}
  \{\dot{v}_v\} =[M_v]^{-1} \left\{ \A_{ \forall p} \Bigl([S_{vp}]^{T}\{b\} V_p - [\nabla_x S_{vp}]^{T}\{\sigma_p\} V_p \Bigr)\right\}.
\end{equation}
where 
\begin{equation}
[M_v] = \A_{ \forall p} \Bigl( [S_{vp}]^{T}[S_{vp}] m_p\Bigl)
\end{equation}
is the consistent mass matrix assembled on the background grid nodes.  The nodal accelerations can be used to increment the material point velocities via a FLuid Implicit Particle (FLIP) update\footnote{FLIP \citep{Brackbill1986} increments the material point velocity field using the change in the vertex velocities over the time step, rather than using the total vertex velocities to overwrite the material point velocity field as with a Particle In Cell (PIC)  \cite{Harlow1964} update.}
\begin{equation}\label{eqn:vp_FLIP}
	\{v_p\}_{n+1} = \{v_p\}_{n} + \Delta t  \sum_{\forall v} S_{vp}\{\dot{v}_v\},
\end{equation}
where $\Delta t$ is the time step size and the subscripts $(\cdot)_n$ and $(\cdot)_{n+1}$ denote the previous and updated states, respectively.  The material point positions are updated using
\begin{equation}\label{eqn:xUpdate}
	\{x_p\}_{n+1} = \{x_p\}_{n} + \Delta t  \sum_{\forall v} S_{vp}\{v_v\}_{n+1}.
\end{equation}
The nodal velocities at the start of the step are obtained by mapping the momentum of the material points to the grid nodes and multiplying this momentum by the inverse of the mass matrix
\begin{equation}\label{eqn:nodeVel}
	\{v_{v}\}_n = [M_v]^{-1}\left\{ \A_{ \forall p} \Bigl(m_p[S_{vp}]^{T}\{v_p\}_n \Bigr)\right\}
\end{equation}
and the nodal velocities at the end of the step are obtained by adding nodal accelerations multiplied by the time step the previous velocity 
\begin{equation}
	\{v_v\}_{n+1}  = \{v_v\}_{n} + \Delta t  \{\dot{v}_v\}.
\end{equation}
The stress updating algorithmic position is an important consideration for explicit material point methods.  There are two options  \cite{Bardenhagen2002}:
\begin{enumerate}
	\item Update Stress First (USF): the stress is updated before the determination of the nodal velocities, (\ref{eqn:nodeAcc}), using the increment in the deformation field from the previous step; or
	\item Update Stress Last (USL): the stress is updated after the determination of the nodal velocities using the increment in the deformation field form the current step. 
\end{enumerate}  
The stress updating point has implications on the energy conservation of MPMs, with the USL approach leading to excessive energy dissipation when combined with a lumped mass matrix \cite{Bardenhagen2002}. However, the USL approach is often preferred as the dissipation is consistent with the accuracy of the solution and it tends to damped unresolved modes \cite{Bardenhagen2002,Berzins2022}.  Differences between USF and USL will be explored in more detail as part of the numerical analyses presented in Section~\ref{sec:num}.

\subsection{Implicit \emph{quasi}-static analysis}\label{sec:IQSA}

For \emph{quasi}-static analysis the inertia term in (\ref{eqn:weak}) is assumed to be negligible, reducing the discrete statement of equilibrium to  
\begin{equation}\label{eqn:MPMweakQS}
  \A_{ \forall p} \Bigl( [\nabla_x S_{vp}]^{T}\{\sigma_p\} V_p - [S_{vp}]^{T}\{b\} V_p \Bigr)= \{0\},
\end{equation}
which is a non-linear statement of equilibrium in terms of the nodal displacements in the deformed configuration.  This non-linear equation can be efficiently solved using an implicit Newton-Raphson approach, which requires the equilibrium statement to be linearised with respect to the nodal displacements to form the global stiffness matrix
\begin{equation}
	[K] =  \A_{ \forall p} \Bigl( [\nabla_x S_{vp}]^{T}[a_p] [\nabla_x S_{vp}]V_p \Bigr),
\end{equation}  
where $[a_p]$ is the spatial consistent tangent modulus of the material point under its current state of deformation (see \citet{Charlton2017} for details).  This global stiffness matrix is used to iteratively update the nodal displacements until the equilibrium statement (\ref{eqn:MPMweakQS}) converges within a given tolerance and at each iteration the material point stress and tangent modulus are updated.  Once convergence has been achieved the material point positions, $\{x_p\}$, are updated using
\begin{equation}
	\{x_p\}_{n+1} = \{x_p\}_{n} +  \sum_{\forall v} S_{vp}\{u_v\}, 
\end{equation}
where $\{u_v\}$ is the displacement of a node over the current load step. This approach is adopted by the open source AMPLE code and the reader is referred to \cite{Coombs2020} for details.

\subsection{Basis functions}

Different material point methods  are characterised by the choice of basis functions, $[S_{vp}]$, that map information between the material points and background grid nodes.  In this paper two options are considered, the:
\begin{enumerate}
	\item standard Material Point Method (sMPM)  \cite{Sulsky1994}, where the shape functions of the underlying finite element grid are adopted and the material points are assumed to be concentrated point masses/volumes; and
	\item Generalised Interpolation Material Point Method (GIMPM)  \cite{Bardenhagen2004}, where the basis functions are formed via the convolution of the shape functions of the finite element grid with a particle characteristic function (usually taken to be a unity function), which generates $C^1$ continuous basis functions via the uniform distribution of the mass/volume associated with the material point over a domain centred on the material point.  In this paper the domains are updated using the symmetric material stretch tensor, as detailed by \citet{Charlton2017}.
\end{enumerate} 
The GIMPM was developed to reduce the well documented cell crossing instability of material point methods, caused by the sudden transfer of internal force as material points cross between elements of the background grid.  Other options include Convected Particle Domain Interpolation (CPDI) methods \cite{Sadeghirad2011,Sadeghirad2013a} and adopting spline-based basis functions \cite{Steffen2008,Andersen2010}.

\subsection{Conditioning issues}\label{sec:test}

\begin{figure}[!h]
\centering
    \begin{subfigure}[t]{0.48\textwidth}
        \centering
        \includegraphics[width=0.98\textwidth]{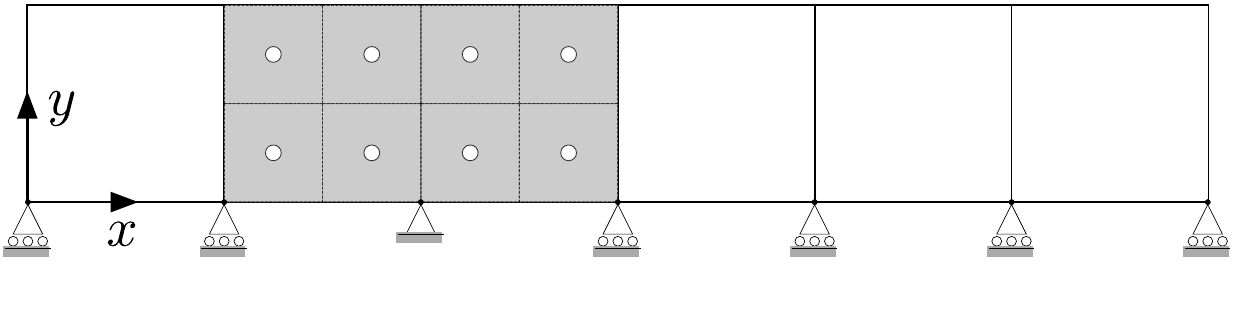}
        \caption{\footnotesize initial position}
    \end{subfigure}%
    \hspace*{0.04\textwidth}
    \begin{subfigure}[t]{0.48\textwidth}
        \centering
        \includegraphics[width=0.98\textwidth]{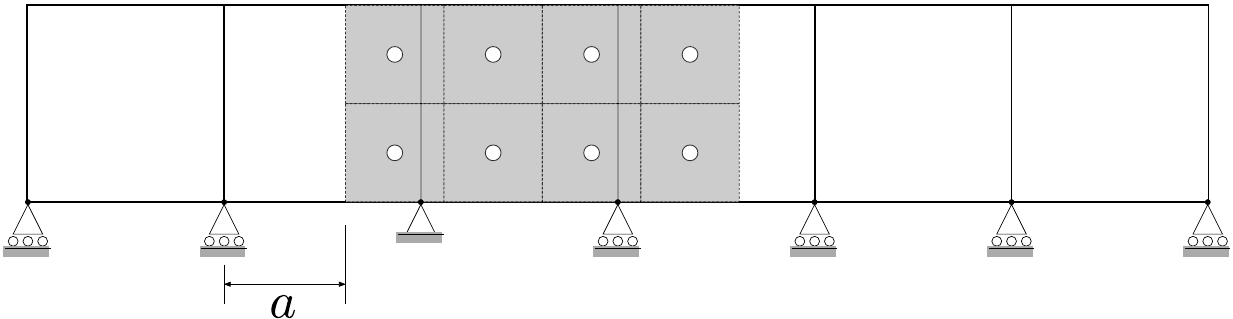}
        \caption{\footnotesize domain translation}
    \end{subfigure}
\caption{Translating material point domain problem setup.}
\label{Fig:testSetup}
\end{figure}

Explicit and implicit material point methods suffer from stability issues associated with the conditioning of the consistent mass and stiffness matrices.  This trivial test problem is designed to highlight these issues.  A $2\times1$m plane strain physical domain was discretised by eight material points, each representing $0.25$m$^3$ of material with a Young's modulus of $E=1$Pa, Poisson's ratio of $\nu=0$ and a density of $\rho=1$kg/m$^3$.  A background grid comprised of $1$m square bi-linear quadrilateral elements was used to construct the global consistent mass and small-strain stiffness matrices.  The background grid nodes were constrained as shown in Figure~\ref{Fig:testSetup}.  

\begin{figure}[!h]
        \centering
    \begin{subfigure}[t]{0.5\textwidth}
        \centering
        \includegraphics[height=62mm]{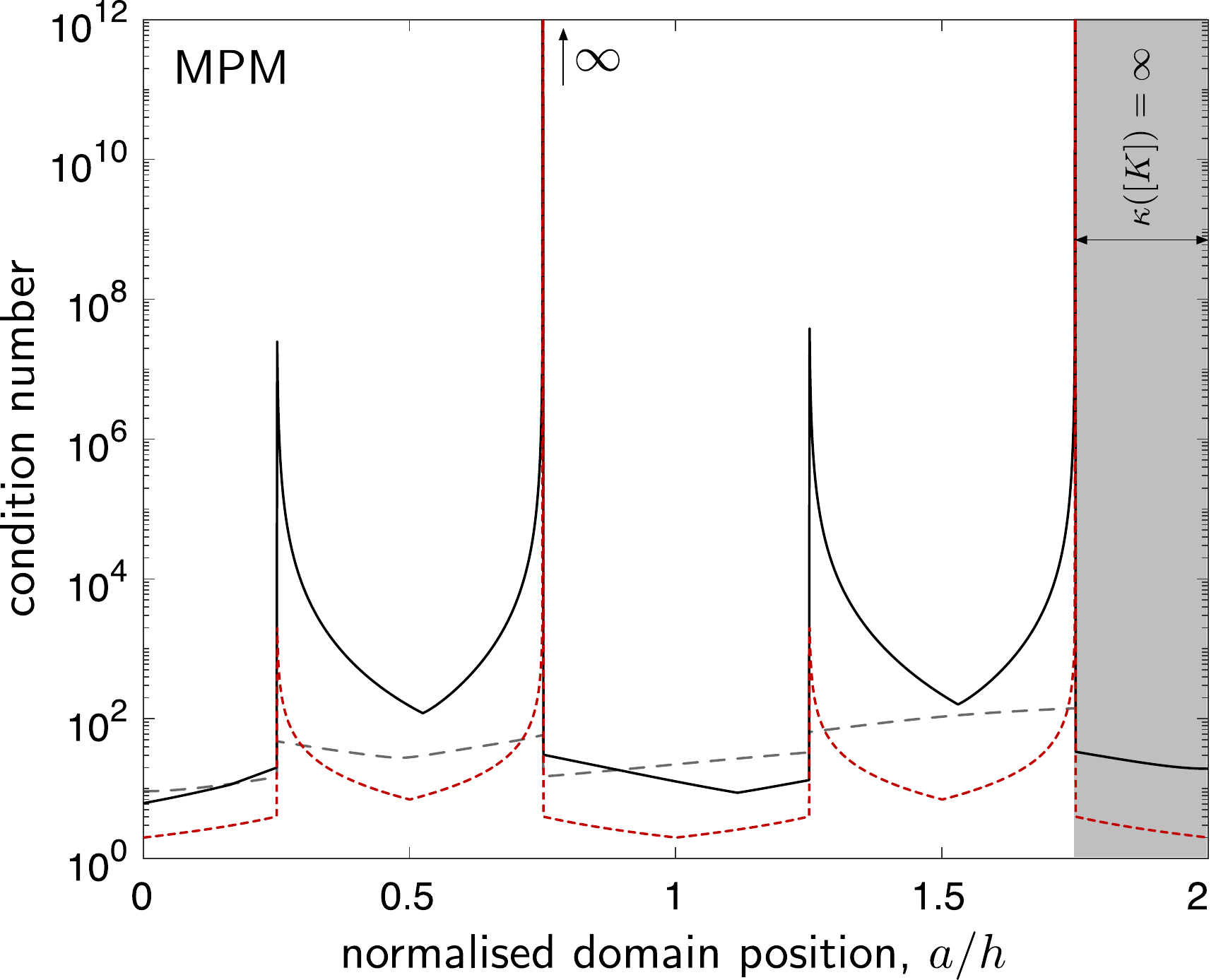}
        \caption{\footnotesize sMPM}
    \end{subfigure}%
    \begin{subfigure}[t]{0.5\textwidth}
        \centering
        \includegraphics[height=62mm]{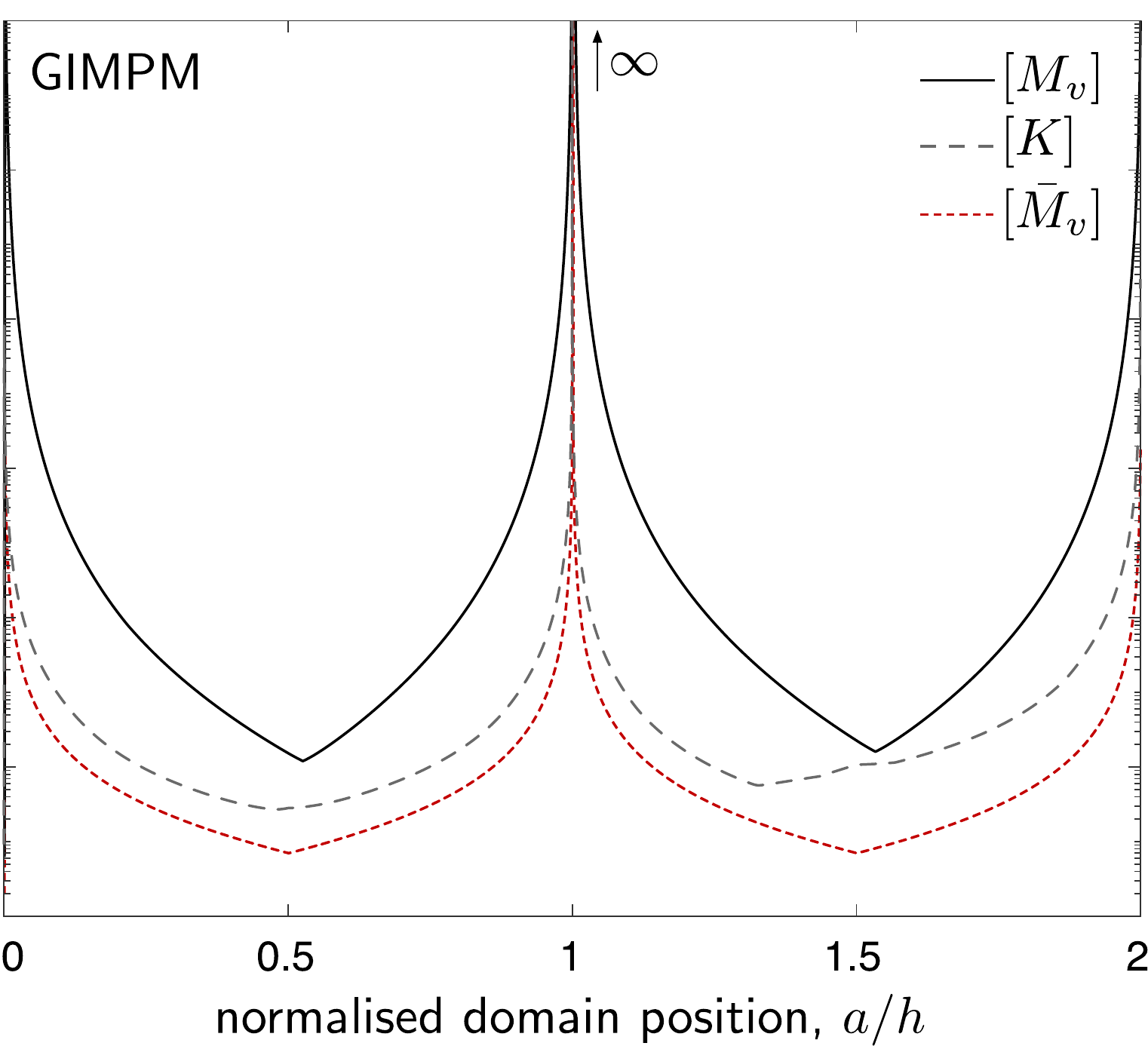}
        \caption{\footnotesize GIMPM}
    \end{subfigure}
	\caption{Translating material point domain: condition number variation.}
	\label{Fig:testCond}
\end{figure}

Figure~\ref{Fig:testCond} shows the evolution of the condition number (the ratio of the largest to smallest eigenvalue\footnote{The analysis presented in this section raises a question regarding what is \emph{reasonable} in terms of the condition number of a linear system of equations.  The answer to this question is related to: the precision of the computational framework being used and the method being used to solve the system of equations.  In all cases, a condition number approaching the reciprocal of the precision of the machine being used will cause accuracy problems.}) of the reduced\footnote{The term \emph{reduced} is referring to a matrix where the rows and columns associated with constrained degrees of freedom have been removed.} consistent mass and stiffness matrices as the physical domain is translated by $a/h=2$ over $5000$ steps for both the standard MPM and GIMPM.  Figure~\ref{Fig:testCond} also reports the condition number of the lumped mass matrix
\begin{equation}
	[\bar{M}_v] = \text{diag}\left( \A_{ \forall p} \Bigl([S_{vp}]^T\{1\} m_p\Bigl) \right)
\end{equation} 
where $\{1\}$ is a vector of ones of length equal to the number of physical dimensions.  The lumped mass matrix is often used in material point methods instead of the consistent mass matrix\footnote{Adopting the lumped mass matrix removes the need to invert the consistent mass matrix in explicit time stepping approaches, replacing it with the inversion of a diagonal matrix, which is trivial.  However, the use of a lumped mass matrix introduces additional numerical dissipation \cite{Burgess1992,Love2006}.  It also does not remove the possibility of spuriously large acceleration values due to very small mass matrix entries.   }.

The condition number of all of the matrices varies as the physical domain translates through the background grid.  For the sMPM (Figure~\ref{Fig:testCond}~(a)), sudden increases in the condition numbers of mass matrices are observed when material points transition between background grid cells (for example at $a/h=0.25,0.75,1.25$ and $1.75$).  This is due to the transitioning material points having very small contributions to certain nodes in terms of the basis functions values.  The stiffness matrix does not show the same spikes as it is formed using the spatial derivatives of the basis functions, which, as we are only concerned with the behaviour in the horizontal direction, are piecewise constant within each element.  Therefore the position of a material point within the element has less impact on the conditioning of the stiffness matrix.  Beyond $a/h=1.75$ the condition of the stiffness matrix $\kappa([K])=\infty$ as there is no constraint in the horizontal direction as all of the material points have moved beyond the influence of the horizontally constrained node.   

The response of the GIMPM is smoother due to the $C^1$ continuity of the basis functions.  However, both the mass and stiffness matrices suffer from conditioning issues as the spatial derivatives of the basis functions are dependent on the overlap between the domains associated with the material points and the background grid cells.  This causes spikes in the condition number when $a/h=0,1$ and $2$.  In all cases the condition number increase is more pronounced for the consistent mass matrix than the lumped mass matrix as it contains the product of the basis function with themselves.  

This is a trivial problem but it demonstrates the conditioning issues faced by material point methods.  The key issue is that we do not know the how the material points, and their associated domains for the GIMPM, will coincide with the background grid through the analysis.  These issues have resulted in most explicit material point methods adopting a lumped mass matrix and/or including thresholds on the minimum considered basis function values that contributes to the analysis.  These numerical \emph{fixes} are not universal and have implications on the accuracy, stability and general use of the methods.  Although there have been some papers investigating the issues associated with poor conditioning of the mass \cite{Ma2010} and stiffness \cite{Wang2016} matrices, a general solution has yet to be provided.

\section{Ghost stabilisation}

\begin{figure}[!h]
        \centering
\begin{subfigure}[t]{0.33\textwidth}
        \centering
        \includegraphics[width=0.99\textwidth]{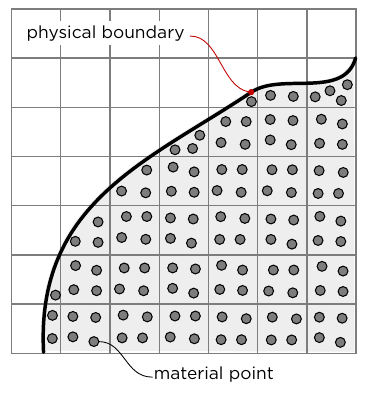}
        \caption{\footnotesize MP discretisation}
    \end{subfigure}%
    \begin{subfigure}[t]{0.33\textwidth}
        \centering
        \includegraphics[width=0.99\textwidth]{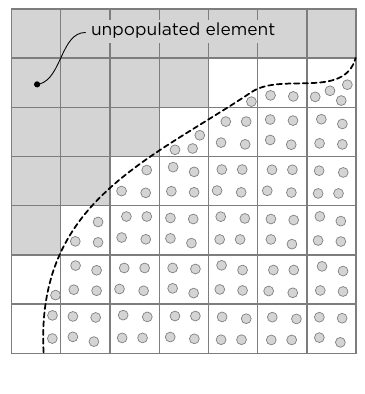}
        \caption{\footnotesize active/inactive elements}
    \end{subfigure}%
    \begin{subfigure}[t]{0.33\textwidth}
        \centering
        \includegraphics[width=0.99\textwidth]{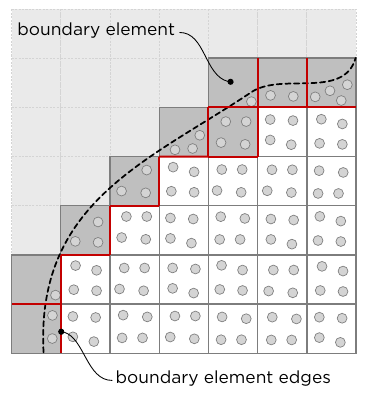}
        \caption{\footnotesize boundary elements}
    \end{subfigure}
	\caption{Material point method discretisation and boundary representation.}
	\label{Fig:boundaries}
\end{figure}

A typical material point situation is shown in Figure~\ref{Fig:boundaries}, where a physical body is discretised by a number of material points on a regular background mesh.  Although a boundary to the physical domain has been shown, most material point methods do not include a representation of the physical boundary and instead they rely on the location of the material points to represent the extent of the physical material.  Figure~\ref{Fig:boundaries}~(b) identifies the unpopulated elements (the elements that do not contain material points) in the background mesh, which have been shaded grey.  Figure~\ref{Fig:boundaries}~(c) identifies the boundary elements (the elements that are intersected by the physical boundary), which are shaded dark grey. It is these boundary elements that potentially contain small overlaps between material points and the background mesh, depending on how the material points are used to represent the physical material and their associated characteristic function.  


One approach to overcome the issues associated with small overlaps between the physical domain and the computational grid is to add additional stabilising terms of the mass/stiffness matrix that provide a bound on the conditioning of the system.    The \emph{Ghost} penalty stabilisation approach was first proposed by \citet{Burman2010}, and aims to bound the conditioning of a matrix by penalising jumps in the gradient of the solution field across elements cut by the physical boundary. The name \emph{Ghost} can be interpreted as projection of the physical solution field into the empty part of the computational mesh.  This is achieved by adding an extra term to the linear system being solved, which strengthens the coupling in the system of equations and ties the solution field in the poorly conditioned part of the linear system to the well conditioned part via the gradient of the solution field near the physical boundary. This extra \emph{Ghost} stabilisation term was expressed by \citet{Sticko2020} (amongst others) as an integral over the faces of elements cut/intersected by the physical boundary, $\Gamma$, with the following bilinear form 
\begin{equation}\label{eqn:stab}
	j(u_i,w_i) = \sum^q_{k=1}  \frac{h^{2k+1}}{(2k+1)(k!)^2} \int_{\Gamma} [[\partial^k_n u_i]] ~[[\partial^k_n w_i]] d\Gamma,
\end{equation}
where $q$ is the polynomial order of the basis functions, $h$ is the size of the element face, $w_i$ and $u_i$ are the test and trial functions ($i$ are the physical dimensions), $\partial^k_n u_i$ is the $k$th spatial derivative of $u_i$ in the direction normal to the face under consideration, $n_i$, and $[[\cdot]]$ denotes the jump over a face, $F$
\begin{equation}
	[[u_i]] = u_i|_{F^+} - u_i|_{F^-}.  
\end{equation}
$F^+$ and $F^-$ denote the faces associated with the \emph{positive} and \emph{negative} elements attached to the boundary between two elements\footnote{The definition of the positive element, $K^+$, and the negative element, $K^-$, associated with a face is arbitrary and swapping the the positive and negative elements will not change of the nature of the stabilisation.  Simply one of the the elements connected to the face is labelled as \emph{positive} and the other as \emph{negative}.}.

\subsection{Boundary identification}

A key aspect of the Ghost stabilisation approach is identifying the element boundaries that are associated with the physical boundary.  The selection of these element boundaries is a clear point of departure between cut-cell finite element and material point methods as most material point simulations do not explicitly define the boundary of the physical domain(s).  Therefore an approach is required to identify these element boundaries without constructing/tracking a representation of the physical boundary.  

\begin{figure}[!h]
        \centering
        \includegraphics[width=0.7\textwidth]{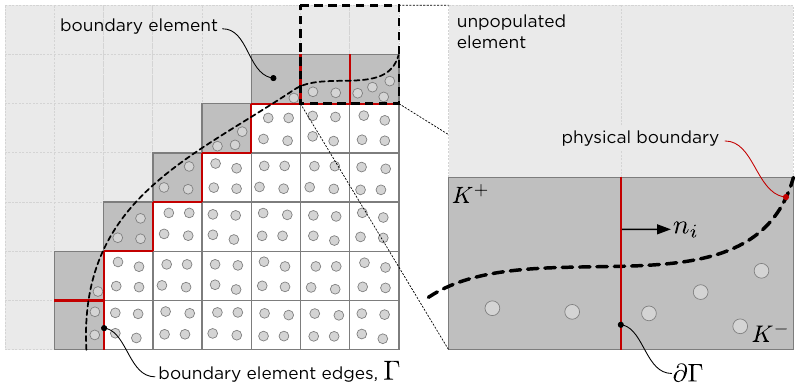}
	\caption{Material Point Method boundaries: positive/negative elements and normal definition.}
	\label{Fig:boundariesDetail}
\end{figure}

In this paper the following steps are taken to identify the element boundaries, $\Gamma$, for the Ghost stabilisation integral:
\begin{enumerate}
	\item Identify the boundary elements, which are the elements \emph{attached} to (sharing a face with) any unpopulated elements (the light grey-shaded elements in Figure~\ref{Fig:boundariesDetail}).  These elements are are the dark grey shaded elements shown in Figure~\ref{Fig:boundariesDetail}.
	\item  The boundary element edges are defined as the boundaries of these elements with: (i) other boundary elements or (ii) elements that are populated by material points (the white-shaded elements in Figure~\ref{Fig:boundariesDetail}).  These element boundaries are shown by the thik red lines in Figure~\ref{Fig:boundariesDetail}\footnote{Note that the dark grey shaded elements and the boundary element edges are different in Figures~\ref{Fig:boundaries}~(c) and \ref{Fig:boundariesDetail} as Figure~\ref{Fig:boundaries}~(c) describes elements \emph{cut} by the boundary rather than the MPM approach described in this section.}.  
\end{enumerate}  
For material point methods it is important to integrate over both the faces between boundary elements and the faces between boundary elements and other \emph{active} elements.  The stabilisation obtained from integrating over the latter faces will impose continuity between the gradient of the solution of the well-conditioned part of the domain and the partially filled elements.  When considering an analysis with multiple bodies, this process should be adopted for each body and then the union of these edges used to determine the stabilisation.

\subsection{Ghost stabilisation for material point analysis}

This paper is restricted to material point methods that adopt bi-linear quadrilateral elements as the underlying background grid.  For these linear elements ($q=1$) the bilinear form of the stabilisation term (\ref{eqn:stab}) becomes
\begin{equation}\label{eqn:stabLinear}
	j(u_i,w_i) = \frac{h^{3}}{3} \int_{\Gamma} [[\partial_n u_i]]~[[\partial_n w_i]] d\Gamma,
\end{equation}
where for a two-dimensional problem the gradient of $u_i$ normal to a boundary can be expressed as 
\begin{equation}
\partial_n u_i = (\nabla u)_{ij} n_j = \left[\begin{array}{cc}
		\frac{\partial u_1}{\partial x} & \frac{\partial u_1}{\partial y} \\ \frac{\partial u_2}{\partial x} & \frac{\partial u_2}{\partial y}
	\end{array} \right]\left\{ \begin{array}{c}
	n_x\\ n_y
	\end{array} \right\}.
\end{equation}
Expanding the jump terms, the stabilisation term for linear elements, (\ref{eqn:stabLinear}), can be expressed as
\begin{equation}
	j(u_i,w_i) = \frac{h^{3}}{3} \int_{\Gamma} \left( \frac{\partial u^+_i}{\partial x_j}n_j- \frac{\partial u^-_i}{\partial x_j}n_j\right) \left( \frac{\partial w^+_i}{\partial x_j}n_j- \frac{\partial w^-_i}{\partial x_j}n_j\right) d\Gamma,
\end{equation}
where $n_j$ is the outward normal for the positive element, as shown in Figure~\ref{Fig:boundariesDetail}. Introducing the finite element approximation space for the test and trial functions gives 
\begin{equation}
	j(u_i,w_i) = \frac{h^{3}}{3} \int_{\Gamma} \left( \frac{\partial [N^+]\{d_u^+\}}{\partial \{x\}}[n]- \frac{\partial  [N^-]\{d_u^-\}}{\partial \{x\}}[n]\right) \left( \frac{\partial [N^+]\{d_w^+\}}{\partial \{x\}}[n]- \frac{\partial [N^-]\{d_w^-\}}{\partial \{x\}}[n]\right) d\Gamma.
\end{equation}
Expanding out this equation and eliminating the nodal values associated with the test function results in a matrix comprised of four sub components multiplied by the physical displacements of the positive, $\{d^+\}$, and negative, $\{d^-\}$, elements 
\begin{equation}
	\{f_G\} = \left[\begin{array}{cc}
		 {[J_G^{++}]}  &  [J_G^{+-}] \\  
		 {[J_G^{-+}]} &  [J_G^{--}]  
	\end{array} \right]\left\{\begin{array}{c}
		\{d^+\}\\ \{d^-\}
	\end{array}\right\} = [J_G]\{d\},
\end{equation}
the combined $[J_G]$ matrix can be compactly expressed as
\begin{equation}\label{eqn:Jmatrix}
	[J_G] = \frac{h^{3}}{3}  \int_{\Gamma} \Bigl([G]^T[m][G]\Bigr) d\Gamma,
	\qquad \text{where} \qquad
	[G] = \bigl[ [G^+] \quad -[G^-] \bigr]^T
\end{equation}
and $[m]=[n][n]^T$.
For two-dimensional analysis, the normal matrix and the matrix containing the derivatives of the basis functions have the following format 
\[
	[n]^T = \left[\begin{array}{cccc}
		n_x & 0 & 0 & n_y \\ 0 & n_y & n_x & 0
	\end{array}\right] 
\qquad \text{and} \qquad
	[G^+] = \left[\begin{array}{cc c cc}
		N^+_1,_x & 0 & \dots & N^+_n,_x & 0\\
		0 & N^+_1,_y  & \dots & 0 & N^+_n,_y \\
		0 & N^+_1,_x  & \dots & 0 & N^+_n,_x \\
		N^+_1,_y & 0 & \dots & N^+_n,_y & 0\\
	\end{array}\right].
\]
$[G^-]$ has the same format as $[G^+]$, with $N^+_i$ replaced with $N^-_i$.  

\subsection{Mass stabilisation for dynamics}

For a dynamic material point method, mass stabilisation can be added to reduce the conditioning issues associated with small overlaps/contributions from material points to the background grid.  This involves adding the following stabilisation matrix to the consistent mass matrix
\begin{equation}\label{eqn:kG}
	[M_G] = \gamma_M [J_G],
\end{equation}
where $\gamma_M$ is a scalar constant that controls how much stabilisation is added to the linear system. \citet{Sticko2020} suggested $\gamma_M = \frac{1}{4}\rho$, where $\rho$ is the density of the material under consideration, is a suitable stabilisation parameter.  The impact of the variation of this parameter will be explored in the numerical analyses in Section~\ref{sec:num}.  Note that the summation of all of the terms in $[J_G]$ is always equal to zero and therefore the addition of the stabilisation does not change the total mass in the system or damage the conservation of momentum of the material point method.  

\subsection{Stiffness stabilisation of quasi-static analysis}\label{sec:stiffStab}

Ghost stabilisation for \emph{quasi}-static analysis is essentially a penalty approach that modifies the \emph{quasi}-static discrete form of the weak form statement of equilibrium (\ref{eqn:MPMweak}) to
\begin{equation}\label{eqn:weakGhost}
    \int_{\varphi_t(K)}[\nabla_x S_{vp}]^{T}\{\sigma_p\} \text{d}V - \int_{\varphi_t(K)}[S_{vp}]^{T}\{b\} \text{d}V  + \beta \int_{\Gamma} [G]^T[n]\{g\}  d\Gamma= \{0\},
\end{equation}
where $\beta=\gamma_k h^3/3$ and $\gamma_k$ is the penalty parameter for Ghost stiffness stabilisation\footnote{The value of this penalty parameter will be discussed later in the section.}, and $\{g\}$ is the jump in the displacement gradient over the boundary element edges, $\Gamma$
\begin{equation}
	\{g\}=[n]^T[G]\{d\}.
\end{equation}
The third term in (\ref{eqn:weakGhost}) is an additional force associated with the Ghost stabilisation that penalises changes in gradient of the solution field over the boundary element edges.  Linearising (\ref{eqn:weakGhost}) with respect to the unknown nodal displacements results in an additional stiffness component linked to this force 
\begin{equation}
	{[K_G]} = \frac{\gamma_k h^3}{3}  \int_{\Gamma} \Bigl([G]^T[m][G]\Bigr) d\Gamma = \gamma_k {[J_G]}.
\end{equation}	
There is significant variation in the value of penalty parameter, $\gamma_k$, used in the literature and also variation in the exponent on the face length, $h$.  For example, when enforcing Dirichlet boundary conditions for elastic wave analysis, \citet{Sticko2020} adopted the following stiffness stabilisation
\begin{equation*}
	[K^D_G] = \frac{\gamma_k}{h^2}[J_G],
\end{equation*}
with the penalty parameter set to $\gamma_k=\mu+\frac{\lambda}{2}$, where  $\lambda$ and $\mu$ are the Lam\'{e} parameters, which sets the stabilisation parameter to be half of the P-wave modulus of the material. \citet{Hansbo2017} suggested a modification to this approach, where a different strength of penalisation was imposed on Neumann (traction) and Dirichlet boundaries.   The Neumann boundary stabilisation was weakened to
\begin{equation*}
 	[K^N_G] = \gamma_k [J_G],
\end{equation*} 
whereas $[K^D_G]$ was adopted on boundaries with Dirichlet (displacement) constraints.  \citet{Hansbo2017} found that $[K^N_G]$ was sufficient to bound the condition number of the linear system, whilst producing more accurate numerical results due to the less onerous penalisation (via the omission of the $h^{-2}$ multiplier on the stiffness stabilisation).  \citet{Hansbo2017} set  $\gamma_k=(2\mu+\lambda)\cdot 10^{-4}$, which is much smaller than that used by \citet{Sticko2020} for elastic wave problems.  \citet{Burman2018} set  $\gamma_k=(\mu+\lambda)\cdot 10^{-7}$ when analysing shape optimisation problems and a number of authors have commented that the results are relatively insensitive to the stabilisation parameter.  This paper is focused on the stabilisation of material point methods due to small overlaps between the background mesh and homogeneous Neumann (traction free) physical boundaries defined by the location of the material points and/or their associated domains and therefore \citet{Hansbo2017}'s Neumann boundary stabilisation is adopted.   The impact of the value of the stabilisation parameter, $\gamma_k$, will be explored in the numerical analysis presented in Section~\ref{sec:num}.

\subsection{Ghost stabilised conditioning}

\begin{figure}[!h]
        \centering
	    \begin{subfigure}[t]{0.5\textwidth}
        \centering
        \includegraphics[height=62mm]{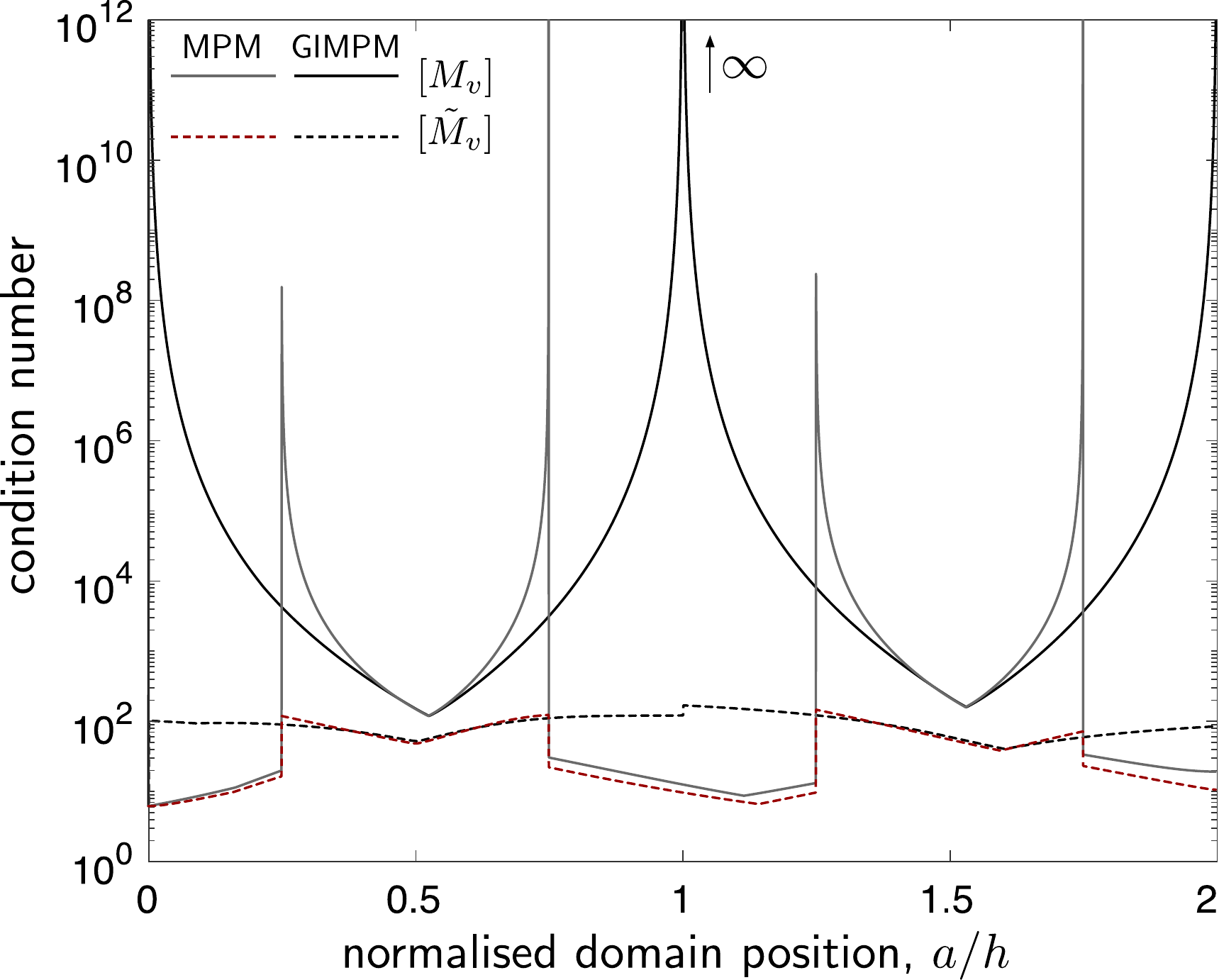}
        \caption{\footnotesize consistent mass matrix}
    \end{subfigure}%
    \begin{subfigure}[t]{0.5\textwidth}
        \centering
        \includegraphics[height=62mm]{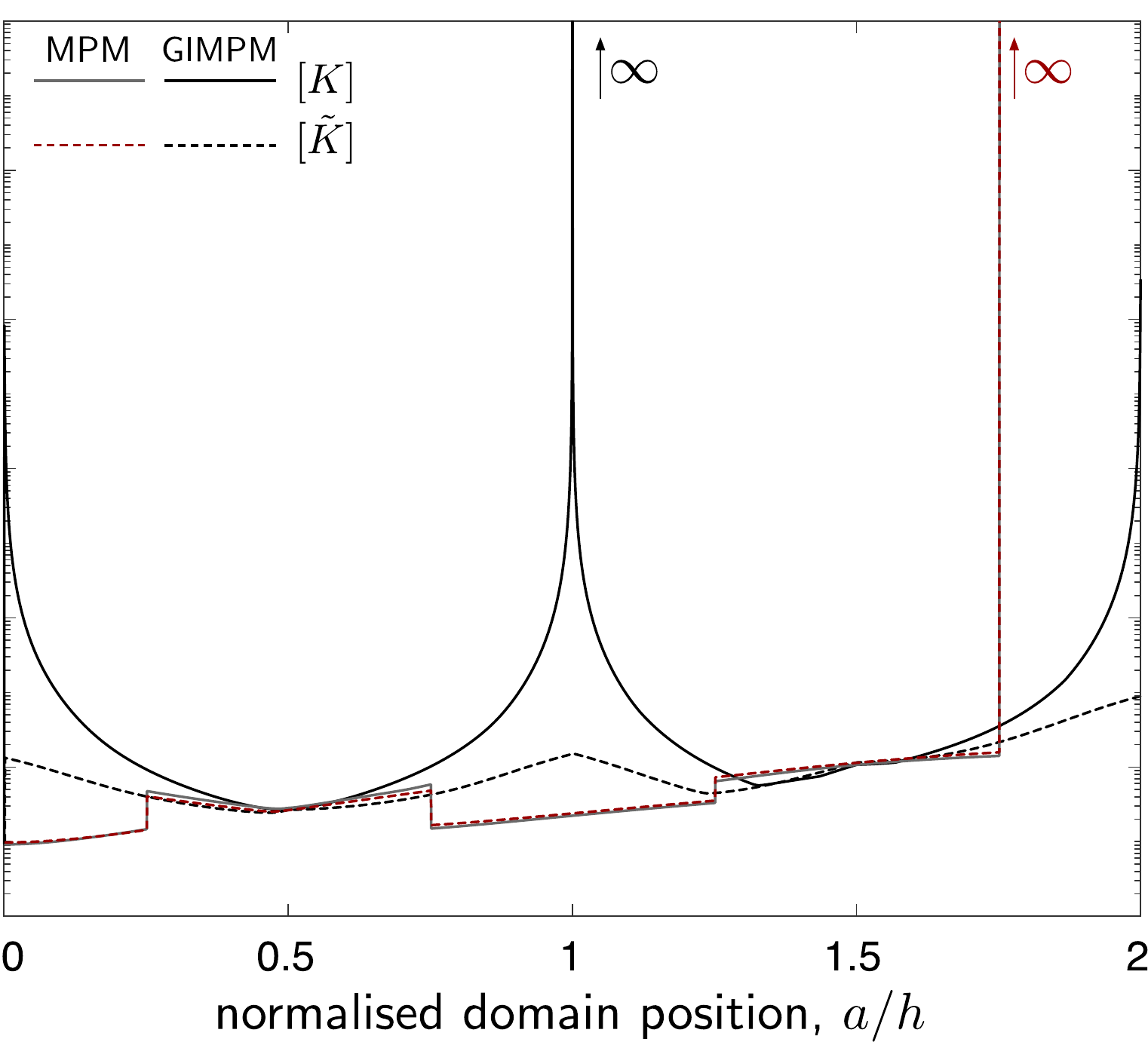}
        \caption{\footnotesize stiffness matrix}
    \end{subfigure}
	\caption{Translating material point domain: condition number variation of (a) the consistent mass and (b) stiffness matrices for the sMPM and GIMPM with and without Ghost stabilisation.}
	\label{Fig:testCondStab}
\end{figure}

The trivial test problem described in Section~\ref{sec:test} is now re-examined with Ghost stabilisation.  Figure~\ref{Fig:testCondStab} shows the variation of the condition numbers of the stabilised, $[\tilde{\cdot}]$, and unstabilised consistent mass and stiffness matrices as the domain translates through the background mesh for the standard MPM and the GIMPM.  The following stabilisation parameters were adopted
\[
	\gamma_M = \frac{1}{4}\rho 
	\qquad \text{and} \qquad 
	\gamma_K=1.
\]  

\begin{table}[!h]
	\centering
	\caption{Translating material point domain: condition number variation of the consistent mass and stiffness matrices for the sMPM and GIMPM with and without Ghost stabilisation. $^{\dag}$maximum for $a/h\in[0,7/4)$.}
	{\small \begin{tabular}{l | c c c c}
		& ~~$\max(\kappa([M_v]))$~~ & ~~$\max(\kappa([\tilde{M}_v]))$~~  & ~~$\max(\kappa([K]))$~~  & ~~$\max(\kappa([\tilde{K}]))$~~\\ \hline
		MPM & $9.08\times10^{26}$ & $1.53\times10^{2}$ & $4.72\times10^{1}$\dag & $4.22\times10^{1}$\dag\\
		GIMPM~~ & $4.43\times10^{43}$ & $1.63\times10^{2}$ & $3.32\times10^{26}$ & $1.66\times10^{2}$ \\
	\end{tabular}}
	\label{Tab:testCondStab}
\end{table}

The Ghost stabilisation removes the spikes in the condition number of the consistent mass matrix for both the standard MPM and the GIMPM.  The is most clearly seen through the maximum condition numbers reported in Table~\ref{Tab:testCondStab}.  The stabilisation has also removed the spikes in the conditioning of the GIMPM stiffness matrix, reducing the maximum condition number by 24 orders of magnitude.  The stiffness matrix of the MPM is well conditioned for $a/h\in[0,7/4)$ and therefore the stabilisation has minimal impact.  The stabilisation does not remove the requirement to constrain rigid body motion and beyond $a/h=7/4$ the condition of the stabilised stiffness matrix is effectively $\kappa([\tilde{K}])=\infty$ as the domain is not constrained in the horizontal direction.

\begin{figure}[!h]
        \centering
        \includegraphics[width=0.6\textwidth]{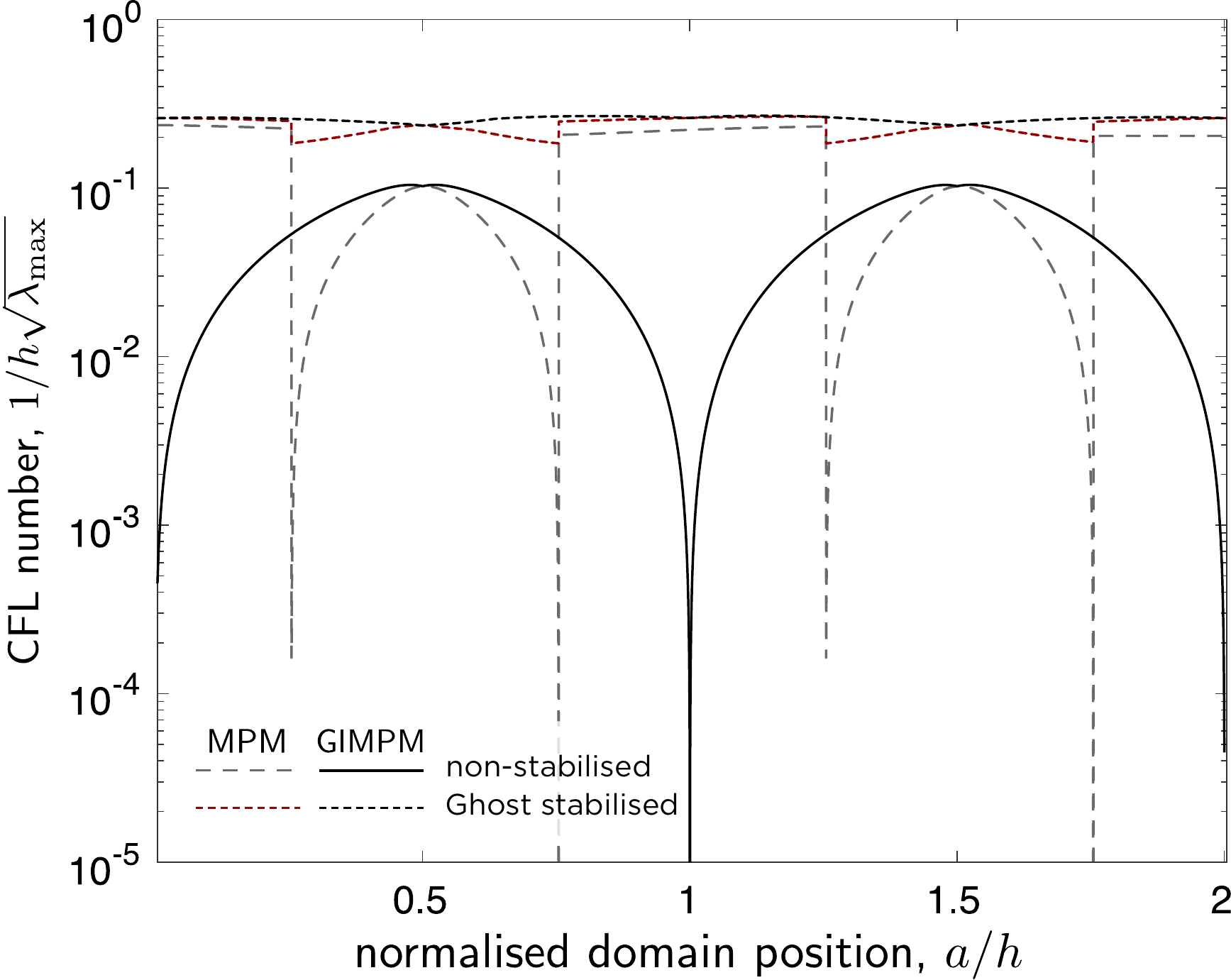}
	\caption{Translating material point domain: CFL number with and without Ghost stabilisation.}
	\label{Fig:testCFL}
\end{figure}

The maximum time step size for an explicit time stepping algorithm is limited by the Courant-Friedrichs-Lewy (CFL) number, $C_{\text{CFL}}$, via 
\[
	\Delta t \leq \alpha C_{\text{CFL}} \min(h),
\]
where $\alpha$ is a constant that depends on the adopted time stepping algorithm \cite{Sticko2020}.  The CFL number can be determined from 
\[
	C_{\text{CFL}} = \frac{1}{h\sqrt{\lambda_{\max}}}, 
\]
where $\lambda_{\max}$ is the largest eigenvalue of the generalised eigenvalue problem $[K]\{x\} = \lambda [M]\{x\}$.  The ideal situation is where the $C_{\text{CFL}}$ number is independent of the position of the physical domain relative to the background mesh and scales linearly with the size of the elements in the background grid.  However, ill-conditioning of the mass and stiffness matrices can cause severe time step limitations  due to the unbounded nature of the matrices' eigenvalues \cite{Sticko2020}.  This issue is demonstrated in Figure~\ref{Fig:testCFL}, where the CFL number is given as the physical domain translates through the background mesh for the sMPM and GIMPM with and without Ghost stabilisation.  The CFL number for the non-stabilised methods are highly dependent on the position of the physical body relative to the background mesh, mirroring the results shown in Figure~\ref{Fig:testCondStab} in terms of the condition numbers of the mass and stiffness matrices.  The CFL number for the Ghost stabilised methods are almost independent of the domain's position on the background grid.  


This conceptual problem has demonstrated the impact of the Ghost stabilisation technique on the conditioning, and therefore accuracy and stability, of the consistent mass and stiffness matrices.  It has taken two matrices that are effectively un-invertible for certain intersections of the background mesh with the physical domain introduced an additional term into the linear system that stabilises the equations and removes the strong dependency of the position of the domain within the background mesh.  It also removes the CFL number dependence on the position physical domain within the background grid.  The impact of this stabilisation will be demonstrated further via the numerical examples presented in Section~\ref{sec:num}, prior to this key aspects of the implementation of the Ghost stabilisation within material point methods will be detailed.

\section{Implementation}

Implementation of the Ghost stabilisation technique only requires minor modification to standard material point method code.  These key differences are described in this section. 

\subsection{Data structures}

The Ghost stabilisation technique requires information on the \emph{skeleton}\footnote{The term \emph{skeleton}, often used in discontinuous Galerkin finite element methods, refers to the boundaries between elements.} of the background mesh that is not routinely stored for the material point method (or continuous Galerkin finite element methods), namely:
\begin{description}
	\item[face connectivity between elements:] this is simply a unique list of all of the internal faces in the analysis and the elements to which they are connected.
	\item[face nodal topology:] although this could be derived from the face connectivity between elements, it is convenient to also store the nodes that are connected to each face, in a similar way to storing the topology of each element in conventional finite element methods.  This face topology information will be used to determine face lengths, normal directions, global quadrature point locations, etc.  
\end{description} 

\subsection{Boundary element edge identification}\label{sec:BEEI}

Determining the boundary element edges is a critical stage of the Ghost stabilisation technique and needs to be performed at each time/load step.  For each physical body:
\begin{enumerate}[(i)]
	\item determine which elements are populated by the material points representing the body\footnote{Determining the active/inactive elements is required for material point methods without stabilisation.}, that is the elements that are \emph{active} in the analysis (the white-shaded elements in Figure~\ref{Fig:boundaries}~(b)/(c));
	\item determine the faces on the boundary between the active and inactive parts of the mesh (faces where one of its associated elements is active and one that is inactive);
	\item using the faces from (ii), identify the \emph{boundary elements} (the dark grey shaded elements in Figure~\ref{Fig:boundaries}~(c));
	\item using the boundary elements from (iii), loop over each face of the element and if both elements associated with the face are \emph{active}, flag as a \emph{boundary element edge} (the thick red lines in Figure~\ref{Fig:boundaries}~(c)).
\end{enumerate} 

\subsection{Ghost stabilisation determination}

The boundary element edges, $\Gamma$, identified using the process described in Section~\ref{sec:BEEI}, must be integrated over to form the Ghost stabilisation matrix (\ref{eqn:Jmatrix}).  Gauss-Legendre quadrature is a convenient and efficient choice, allowing $[J_G]$ to be approximated as
\begin{equation}
	[J_G] \approx  \A_{ \forall \Gamma} \left( \frac{h^{3}}{3}  \sum_{i=1}^{n_{Gp}} \Bigl([G_i]^T[m][G_i]\Bigr) w_i \det([\mathcal{J
}_i]) \right),
\end{equation}
where the $(\cdot)_i$ subscripts denote quantities that are potentially dependent on the Gauss point location, $w_i$ are the Gauss point weights, $\det([\mathcal{J}])$ is the determinant of the surface Jacobian (the ratio of the global/local face lengths, $h/2$ for two-dimensional linear elements) and $n_{Gp}$ is the number of Gauss points used to approximate the integral.

The implicit \emph{quasi}-static approach described in Section~\ref{sec:IQSA} adopts an updated Lagrangian approach and therefore both the normal to each face and the spatial derivatives of the basis functions will change over the load step.  However, there would be minimal benefit to including this variation as the Ghost stabilisation is not a physical quantity.  The \emph{cost} of including this variation would be that the stabilisation matrix and the associated force contribution to the weak statement of equilibrium would need to be recalculated at each Newton-Raphson iteration within every load step.  Assuming that the additional stiffness is constant over a load step also allows the force contribution from the Ghost stabilisation to the weak statement of equilibrium to be calculated via the product of the Ghost stabilisation stiffness matrix, $[K_G]$, with the background grid nodal displacements from the current step.  In summary, as the stabilisation is a numerical fix to improve the stability of the method, it is assumed that $[J_G]$ is constant over each load/time step.

\section{Numerical investigations}\label{sec:num}

All of the simulations presented in this section are conducted under the assumption of plane strain and adopt a background grid of bi-linear quadrilateral elements.  The generalised interpolation material point method is adopted for all analyses.  

\subsection{Mass stabilisation: translating, rotating and stretching}

The first set of numerical analyses are focused on mass stabilisation and key steps that are required in dynamic material point simulations, namely the mapping of velocity between the material points and the grid at the start of a time step. 

\subsubsection{Mass stabilisation: rigid body translation}\label{sec:RBT}

\begin{figure}[!h]
        \centering
	    \begin{subfigure}[t]{0.4\textwidth}
        \centering
        \includegraphics[width=\textwidth]{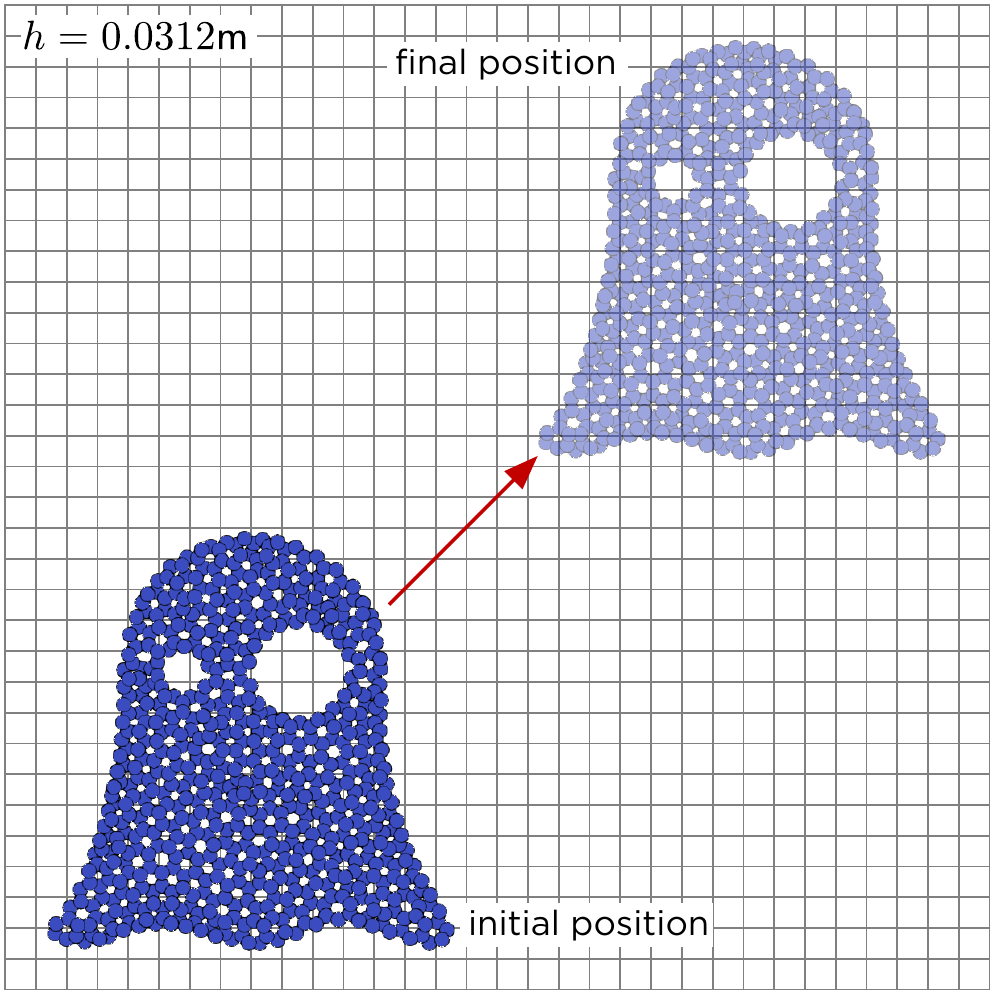}
        \caption{\footnotesize}
    \end{subfigure}%
    \begin{subfigure}[t]{0.4\textwidth}
        \centering
        \includegraphics[width=\textwidth]{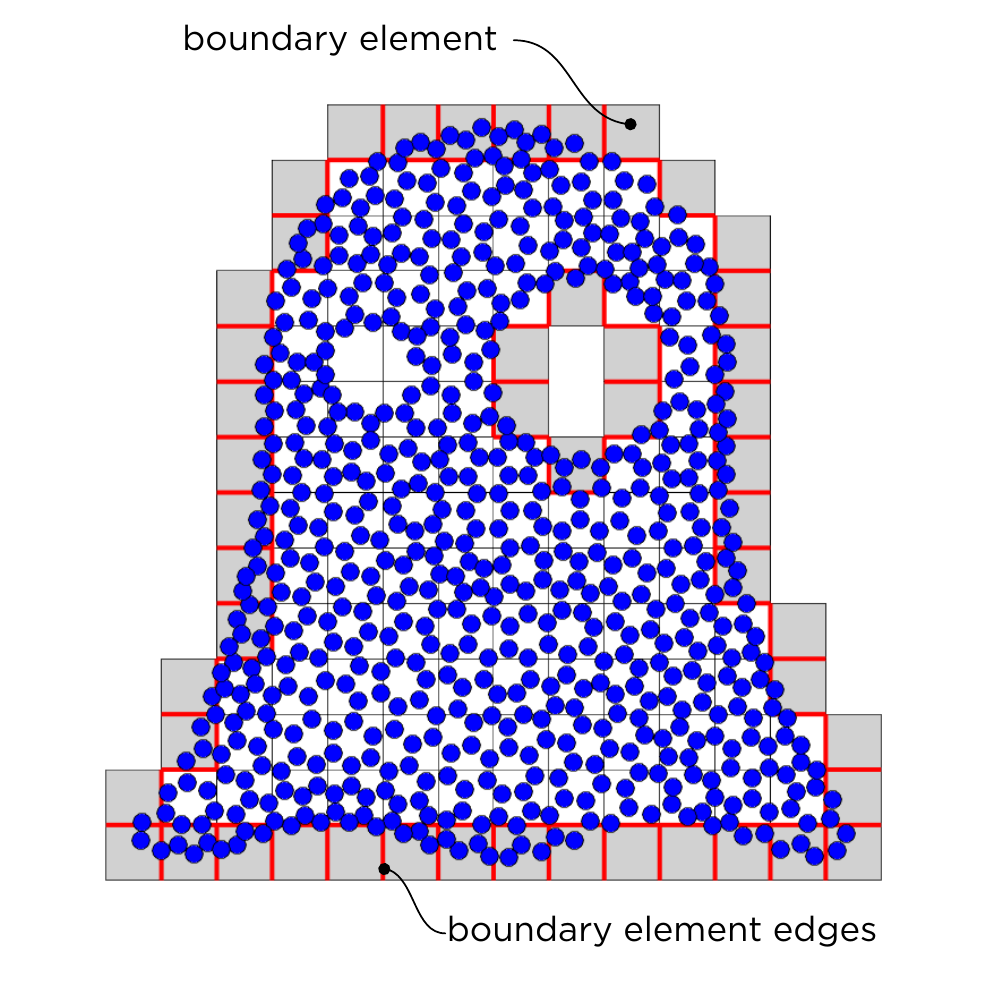}
        \caption{\footnotesize }
    \end{subfigure}
	\caption{Translating ghost discretisation: (a) initial and final material point positions and (b) boundary elements (grey shaded) and the faces (thick red lines) for the initial material point distribution (only the \emph{active} elements in the background mesh are shown).}
	\label{Fig:Ghost}
\end{figure}

This experiment investigates the variation of the condition number of the consistent mass matrix of the active elements as a ghost emoji translates through a regular background mesh.  A unit background mesh of bi-linear square elements was adopted and the ghost had a height and width of $0.4$m, as shown in Figure~\ref{Fig:Ghost}.  The ghost underwent a rigid body motion of 
\[
	\{u\} = \{0.5 \quad 0.5\}^T~\text{m}
\]  
over 1000 steps and the condition number of the global (non-stabilised) consistent mass matrix, $[M_v]$, was recorded at each step along with the condition of the stabilised, $[\tilde{M}_v]$, and lumped, $[\bar{M}_v]$, mass matrices.  The Ghost had a density of $\rho=1$kg/m$^3$ and was represented by 646 generalised interpolation material points\footnote{The positions, volume, mass, etc. of all of the material points for all of the analysis presented in the paper are provided as VTK files within the supplementary data associated with the paper.  See the Acknowledgements section for details.}, as shown in Figure~\ref{Fig:Ghost}.  

\begin{figure}[!h]
\centering
    \begin{subfigure}[t]{0.5\textwidth}
        \centering
        \includegraphics[width=0.98\textwidth]{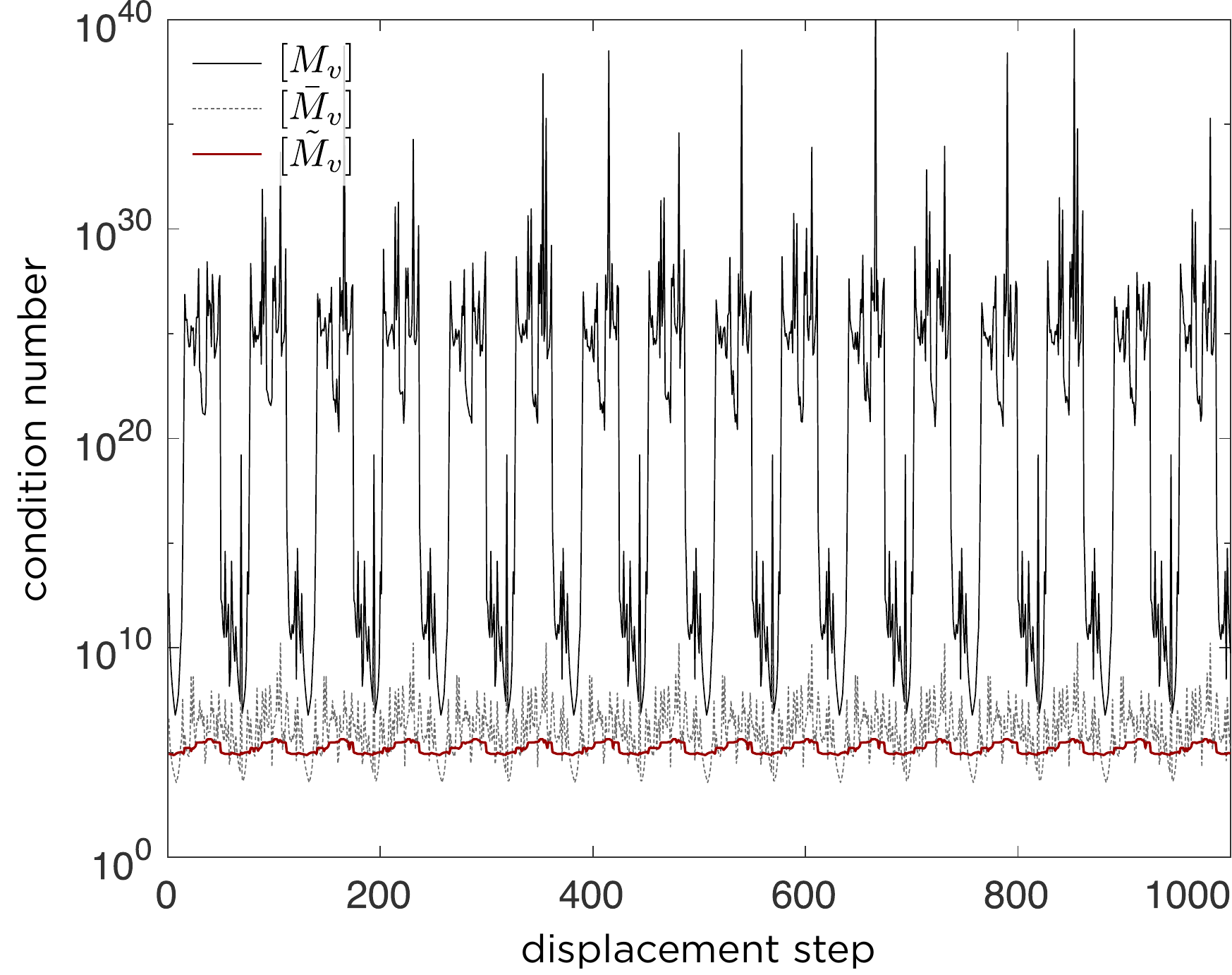}
    \end{subfigure}%
    \begin{subfigure}[t]{0.5\textwidth}
        \centering
        \includegraphics[width=0.98\textwidth]{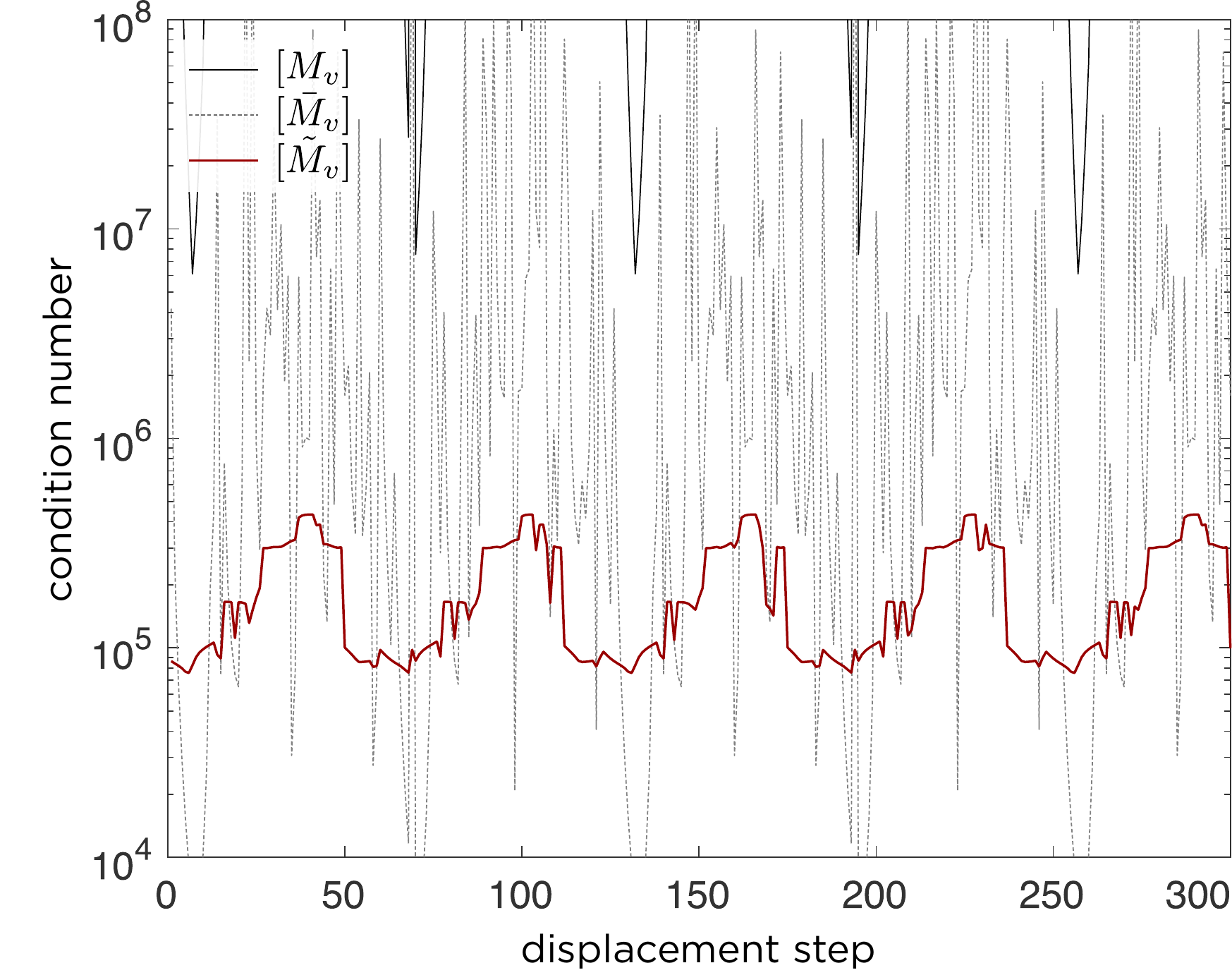}
    \end{subfigure}
\caption{Translating ghost: mass matrix condition number variation with displacement step.}
\label{Fig:GhostCond}
\end{figure}

Figure~\ref{Fig:GhostCond} shows the evolution of the condition number of the consistent (black line), stabilised (thick red line) and lumped (grey dashed line) mass matrices as the Ghost translates through the background mesh with $h=0.0312$m and a mass stabilisation parameter of $\gamma_M=\frac{1}{4}\rho$.  The maximum condition number of the consistent mass matrix was $1.87\times10^{41}$, whereas the stabilised mass matrix had a maximum value of $4.33\times10^{5}$.


%


\begin{table}[!h]
\centering
\caption{Translating ghost: maximum mass matrix condition number for different numbers of material points per Ghost with $h=0.0312$ and $\gamma_M=0.25\rho$.}
{\small \begin{tabular}{l | ccc}
						& ~~~~~coarse~~~~~ & ~~~~~medium~~~~~ & ~~~~~~fine~~~~~~ \\ 
	no. material points/Ghost~ 		& 646 				& 2,389 				& 9,613 \\ \hline 
	non-stabilised $[M_v]$ 		& $1.87\times10^{41}$	& $2.38\times10^{36}$	& $3.01\times10^{33}$\\
	stabilised $[\tilde{M}_v]$ 		& $4.33\times10^{5}$		& $3.12\times10^{5}$ 		& $2.97\times10^{5}$\\
	lumped (diagonal) $[\bar{M}_v]$ 	& $1.71\times10^{10}$	& $7.66\times10^{11}$	& $1.20\times10^{11}$\\ \hline
\end{tabular}}
\label{Tab:Ghost_noMPs}
\end{table}

Table~\ref{Tab:Ghost_noMPs} explores the variation of the maximum condition number of the non-stabilised, stabilised and lumped mass matrices with different numbers of generalised interpolation material points with $\gamma_M=0.25$.  Due to the significant difference between the stabilised and non-stabilised matrices, the condition number of the stabilised matrix scales approximately linearly with the stabilisation parameter, $\gamma_M$.  That is, reducing/increasing the stabilisation parameter by an order of magnitude will increase/reduce the condition number of the stabilised mass matrix by an order of magnitude until the condition number of $[J_G]$ approaches $[M_v]$.


\subsubsection*{Velocity mapping}

\begin{figure}[!h]
\centering
\includegraphics[width=0.5\textwidth]{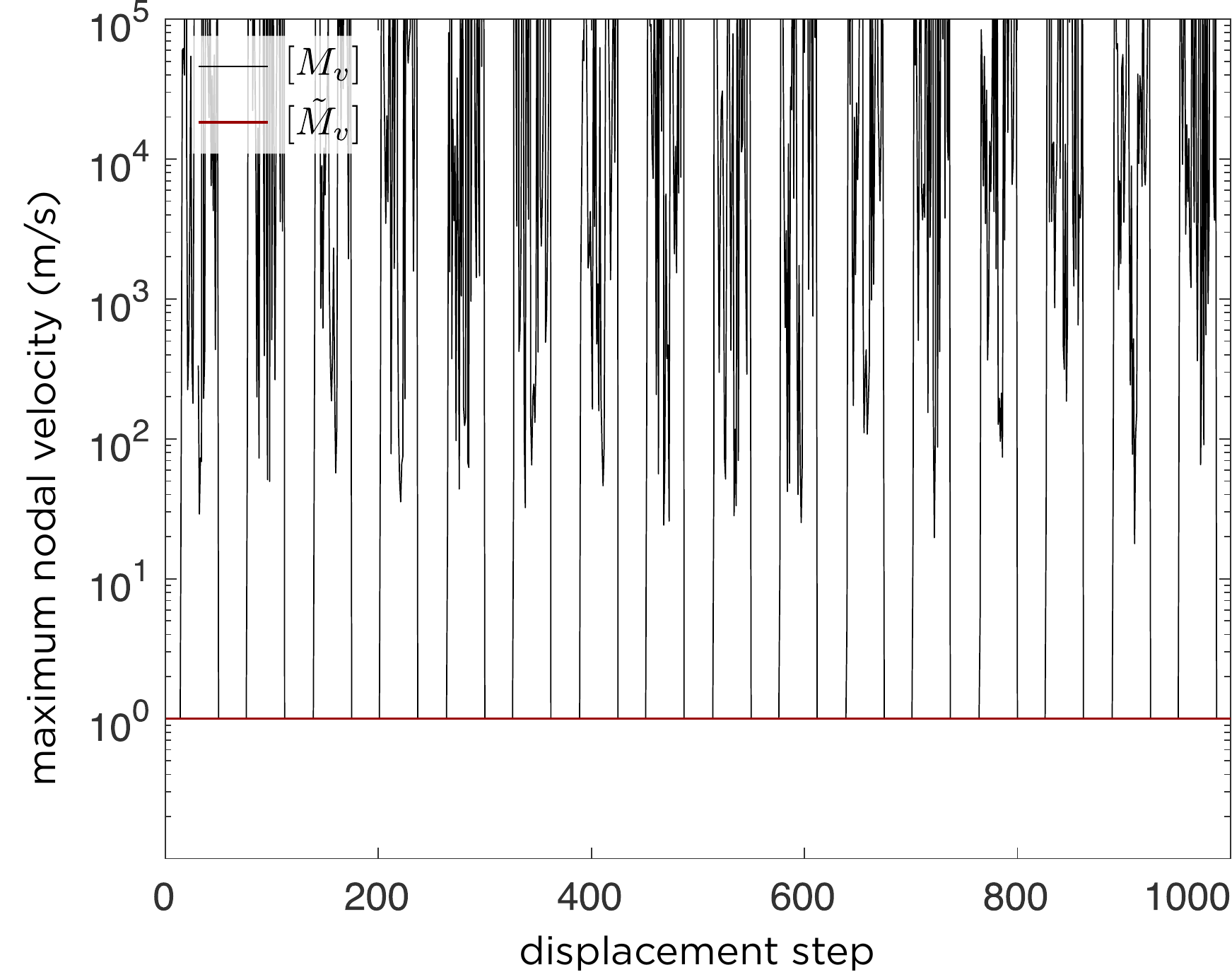}
\caption{Translating ghost: maximum nodal velocity variation with displacement step.}
\label{Fig:GhostVel}
\end{figure}

The stabilised and non-stabilised consistent mass matrices were also used to map the velocity from the coarse material point distribution to nodes using (\ref{eqn:nodeVel}).  In this numerical experiment $\{v_p\}$ was set to
\[
	\{v_p\} = \{1 \quad 1\}^T
\]
and therefore to be consistent with this constant velocity field, all \emph{active} nodes in the background mesh should have a velocity equal to $\{v_v\}=\{1 \quad 1\}^T$.  The variation of the maximum nodal velocity with displacement step for the stabilised (red line) and non-stabilised consistent (black line) mass matrices are shown in Figure~\ref{Fig:GhostVel}.  The nodal values of velocity obtained from the non-stabilised consistent mass matrix are highly dependent on the position of the ghost relative to the background mesh, whereas the stabilised mass matrix correctly predicts nodal velocities of  $\{1 \quad 1\}^T$, with a maximum error of $3.32\times10^{-13}$, irrespective of the position of the ghost within the background mesh.  The lumped mass matrix correctly predicts the nodal velocities as, for constant velocity fields, the summation of the mass basis functions in the lumped matrix and nodal momentum cancel out irrespective of the basis function values, allowing the lumped mass matrix to achieve machine precision in the nodal velocities with a maximum error of $1.11\times10^{-15}$.  This is however only true for velocity fields that are constant in space, as will be demonstrated in the next section.

\subsubsection{Explicit dynamics: stretching ghost \& velocity update}

\begin{figure}[!h]
\centering
    \begin{subfigure}[t]{0.35\textwidth}
        \centering
        \includegraphics[width=\textwidth]{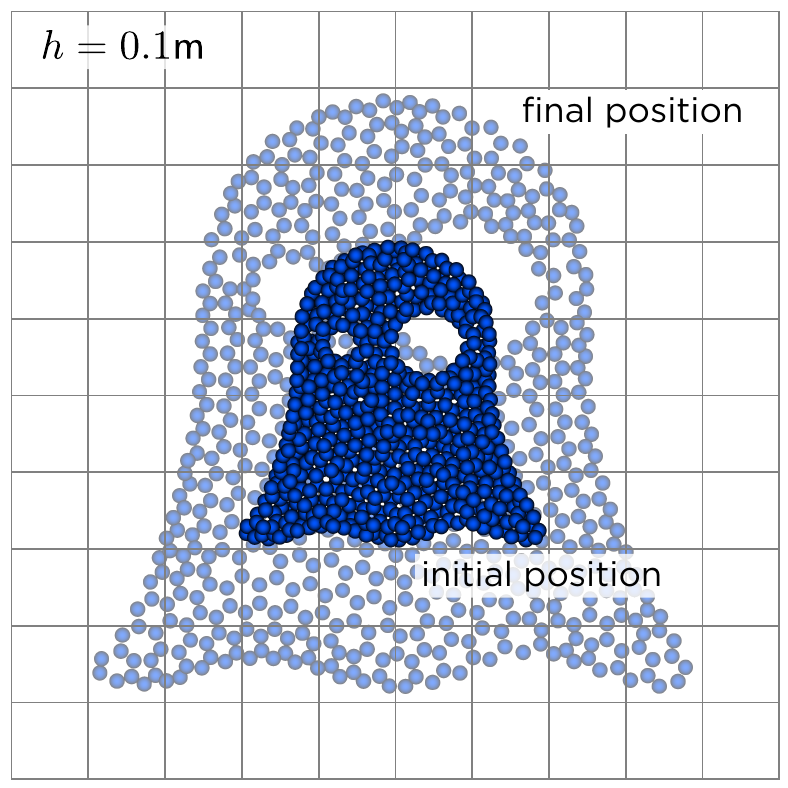}
        \caption{\footnotesize discretisation}
    \end{subfigure}\hspace*{5mm}%
    \begin{subfigure}[t]{0.35\textwidth}
        \centering
        \includegraphics[width=\textwidth]{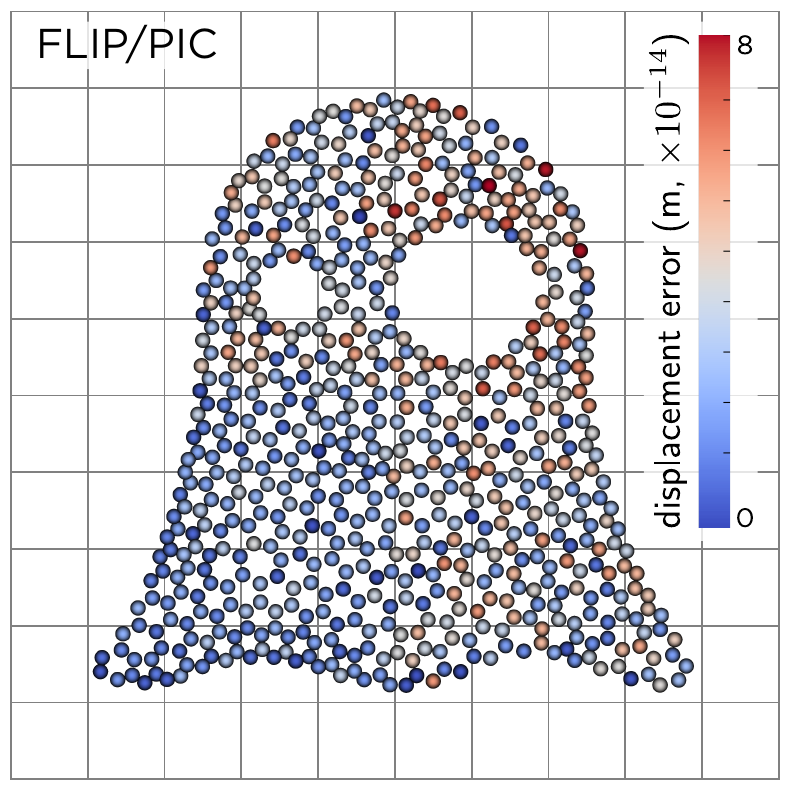}
        \caption{\footnotesize ghost stabilised}
    \end{subfigure}\\[2mm]
    \begin{subfigure}[t]{0.35\textwidth}
        \centering
        \includegraphics[width=\textwidth]{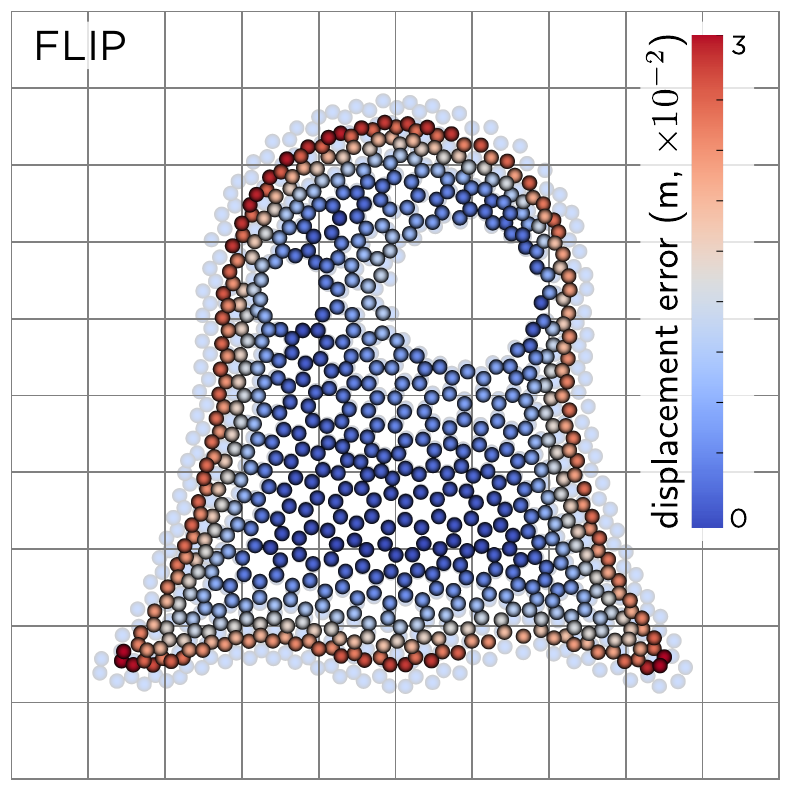}
        \caption{\footnotesize lumped FLIP}
    \end{subfigure}\hspace*{5mm}%
     \begin{subfigure}[t]{0.35\textwidth}
        \centering
        \includegraphics[width=\textwidth]{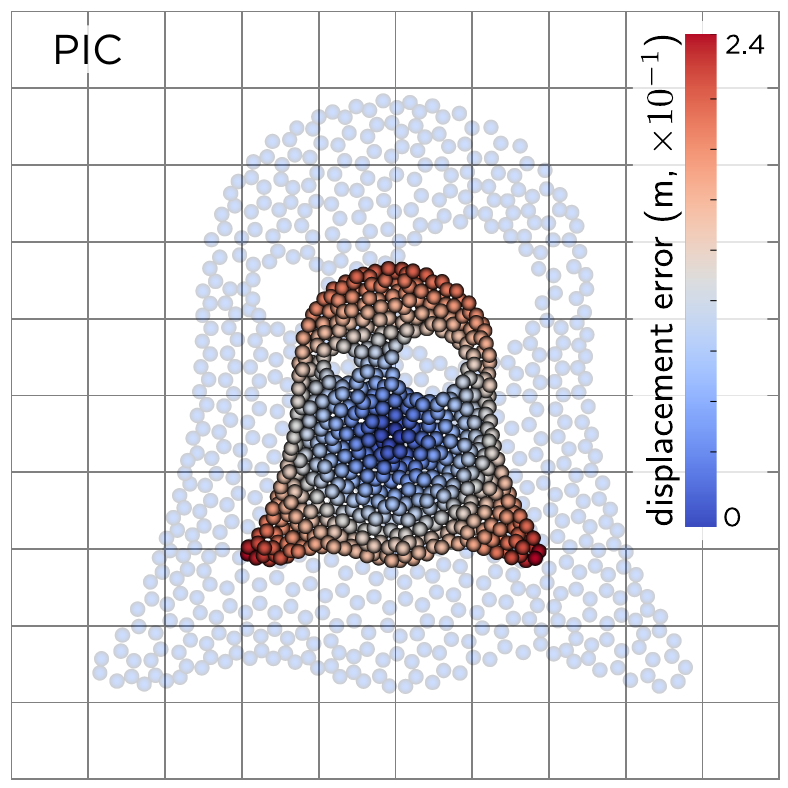}
        \caption{\footnotesize lumped PIC}
    \end{subfigure}
\caption{Stretching ghost: (a) initial discretisation and expected final material point positions for $h=0.1$m and the material point positions, coloured according to the norm of the displacement error, at the end of the analysis for: (b) the Ghost stabilised consistent mass matrix and the lumped mass matrix using a (c) FLIP and (d) PIC material point velocity update.}
\label{Fig:StretchGhost}
\end{figure}

This experiment investigates the differences in the velocity mapping process using the consistent and lumped mass matrices as a ghost emoji expands on a regular background mesh using the explicit dynamic material point formulation described in Section~\ref{sec:ExpDyn}.  A unit background mesh  with $h=0.1$m was adopted.  The ghost had a height and width of $0.4$m and was centred within the background mesh, as shown by the dark blue points in Figure~\ref{Fig:StretchGhost}~(a).  The \emph{coarse} material point discretisation (646 material points), as described in the previous section, was adopted for this analysis. For this artificial problem, the ghost emoji had zero stiffness ($E=0$Pa) and a density of $\rho=1000$kg/m$^3$. The effect of gravity was ignored and the behaviour modelled over 1 second using 1,000 time steps.  A mass stabilisation parameter of $\gamma_M=\frac{1}{4}\rho$ was adopted for the Ghost-stabilised simulations.
The problem was initialised such that the material points had a velocity of
\[
	\{v^0_p\} = \frac{1}{t_t}\Big\{ \{x_p\} - \{x_{\text{cen}}\} \Bigr\},
\]
where $\{x_{\text{cen}}\}=\{0.5 \quad 0.5\}^{T}$m was the centre of ghost emoji and $t_t=1$s the total simulation time.  This velocity field represents a uniform expansion about $\{x_{\text{cen}}\}$.   As the ghost has zero stiffness, there is nothing to resist the expansion of the body and the material points should move with a constant velocity throughout the analysis, doubling the size of the ghost over 1 second (as shown by the light blue points in Figure~\ref{Fig:StretchGhost}~(a)).  This means that at any point in the analysis, the norm of the displacement error for any material point in the body can be evaluated using the following Euclidean norm
\[
	u^p_{\text{error}} = \Bigl|\{u_p\} - t\{v_p^0\}\Bigr|_2,
\] 
where $\{u_p\}$ is the displacement of a material point at time $t$.  In terms of deformation, all material points should have a deformation gradient of $[F]=(1+t)[I]$ with a volume ratio of $J=(1+t)^2$.  Note that it was not possible to use the non-stabilised consistent mass matrix for this analysis, with the simulation diverging after 21 time steps.  

Figures~\ref{Fig:StretchGhost}~(b) and (c) show the material point positions, coloured according to the norm of the displacement error, at the end of the analysis for the (b) Ghost stabilised consistent and (c) lumped mass matrices using a FLIP velocity update (\ref{eqn:vp_FLIP}).  It is worth highlighting that the scales on the colour bars are different for the two figures, with Figure~\ref{Fig:StretchGhost}~(c) having a range of $[0,3]\times10^{-2}$ whereas Figure~\ref{Fig:StretchGhost}~(b) has a range 12 orders of magnitude smaller, at $[0,8]\times10^{-14}$.  Figure~\ref{Fig:StretchGhost}~(c) also shows the expected positions of the material points at the end of the analysis using light blue shaded circles.    Figure~\ref{Fig:StretchGhost}~(d) shows the final material point positions using a PIC material point velocity update\footnote{The PIC approach \cite{Harlow1964} uses the total nodal velocities to overwrite the material point velocity field, assuming that the velocity through the domain varies according to the vertex values and their associated basis functions.}, where again the material points are coloured according to the norm of the displacement error.  A number of points can be observed from the results presented in Figures~\ref{Fig:StretchGhost}~(b)-(d):
\begin{enumerate}[(i)]
	\item for this zero acceleration problem, a FLIP material point velocity update maintains the correct velocities at the material points throughout the analysis irrespective if a lumped of consistent mass matrix is adopted \cite{Pretti2022}; 
	\item the Ghost stabilised consistent mass matrix correctly predicts the material point velocities throughout the analysis with a displacement error of less than $10^{-13}$m across all material points and time steps;
	\item the lumped mass matrix causes in an error in the updated material point positions for non-rigid body motions due to errors in the mapping of the material point velocities to the background grid vertices using (\ref{eqn:nodeVel}) as these nodal velocities are used to update the material point positions at the end of the time step such that the deformation of the material points are consistent with the background grid.   Therefore adopting a lumped mass matrix means that the material point velocities are no longer consistent with their deformation over the time step, which is mapped from the grid incremental displacements via (\ref{eqn:xUpdate});
	\item the Ghost stabilised consistent mass matrix predicts vertex velocities, $\{v_v\}$, that are consistent with the physical linear velocity variation due to the Ghost stabilisation penalising variations of the gradient of the basis functions across the boundary element edges this means that the method predicts the correct material point displacements (and velocities) using both FLIP and PIC velocity updates; and
	\item combining a PIC velocity update with a lumped mass matrix destroys the initial material point velocity field due to inaccuracies in the vertex velocities (as explained under point (ii) and also by \citet{Pretti2022}), leading to a physically unrealistic solution as the total velocity field is repeatedly mapped between the material points and the grid vertices.
\end{enumerate}
This simple problem has demonstrated the benefits of adopting a Ghost stabilised consistent mass matrix over a lumped mass matrix.  It is important to highlight that adopting a lumped mass matrix also has consequences for the energy conservation of material point methods (as discussed in detail by \citet{Love2006}).  This point is explored in the following numerical example. 

\subsection{Explicit dynamics: collision of elastic ghosts}

\begin{figure}[!h]
	\centering
	\includegraphics[width=0.35\textwidth]{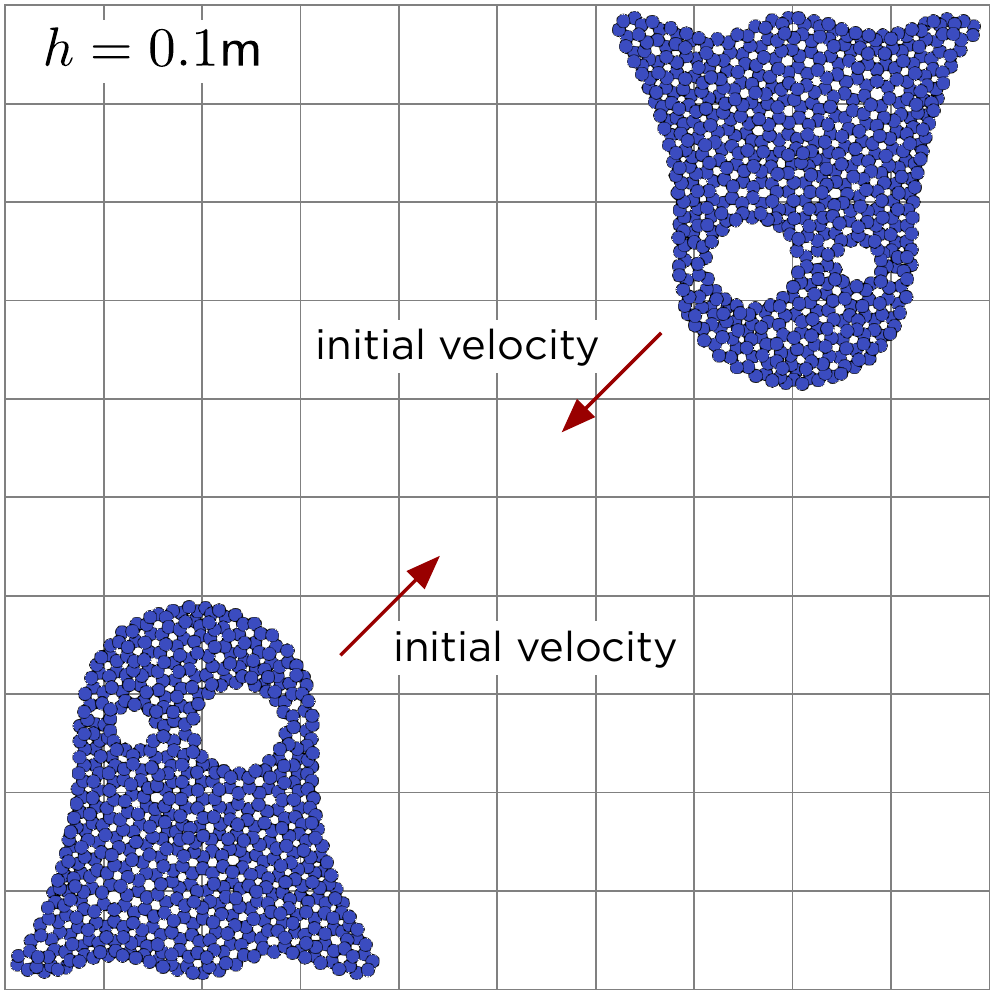}
	\caption{Collision of elastic ghosts:  initial \emph{coarse} material point  discretisation with $h=0.1$m.}
	\label{Fig:ghostCollisionSetup}
\end{figure}

This example considers the collision of two elastic ghosts and has been designed to investigate the energy conservation implications of adopting a lumped or a stabilised consistent mass matrix.  The analyses were run on a 1m$\times$1m background grid, with 5, 10 and 20 elements in each direction ($h=0.20$m, $0.10$m and $0.05$m). The ghosts were initially placed in the upper right and lower left corners of the background grid, with velocities of $\{-0.1,-0.1\}$m/s and $\{0.1, 0.1\}$m/s respectively (as shown in Figure~\ref{Fig:ghostCollisionSetup}). Each ghost had width and height of $0.4$m, a Young's modulus of $E=1000$Pa, a Poisson's ratio of $\nu=0.3$ and a density of $\rho=1000$kg/m\textsuperscript{3}.  The same three generalised interpolation material point resolutions used in the previous sections, and detailed in Table~\ref{Tab:Ghost_noMPs}, were adopted for this analysis.  The effect of gravity was ignored and the behaviour modelled over 3.5 seconds using a total of $n_t=1000$, $2000$ and $4000$ time steps. The system was undamped and was analysed using the explicit dynamic material point formulation described in Section~\ref{sec:ExpDyn} with a lumped and stabilised consistent mass matrix ($\gamma_M=\frac{1}{4}\rho$). Both Update Stress First (USF) and Update Stress Last (USL) material point stress updating procedures were considered.  The total initial energy in the system for all of the analyses was $W_0=2.14726$J.  

Note that it is not possible to analyse this problem using a non-stabilised consistent mass matrix for any of the combinations of parameters used with the lumped/stabilised analyses.  When trying to use the non-stabilised consistent mass matrix, the simulation becomes unstable after a small number of time steps and the ghosts \emph{explode} due to spurious nodal velocity values associated with the inversion of $[M_v]$.

\begin{figure}[!h]
\centering

\begin{subfigure}[c]{0.49\textwidth}
        \centering
        \includegraphics[width=\textwidth]{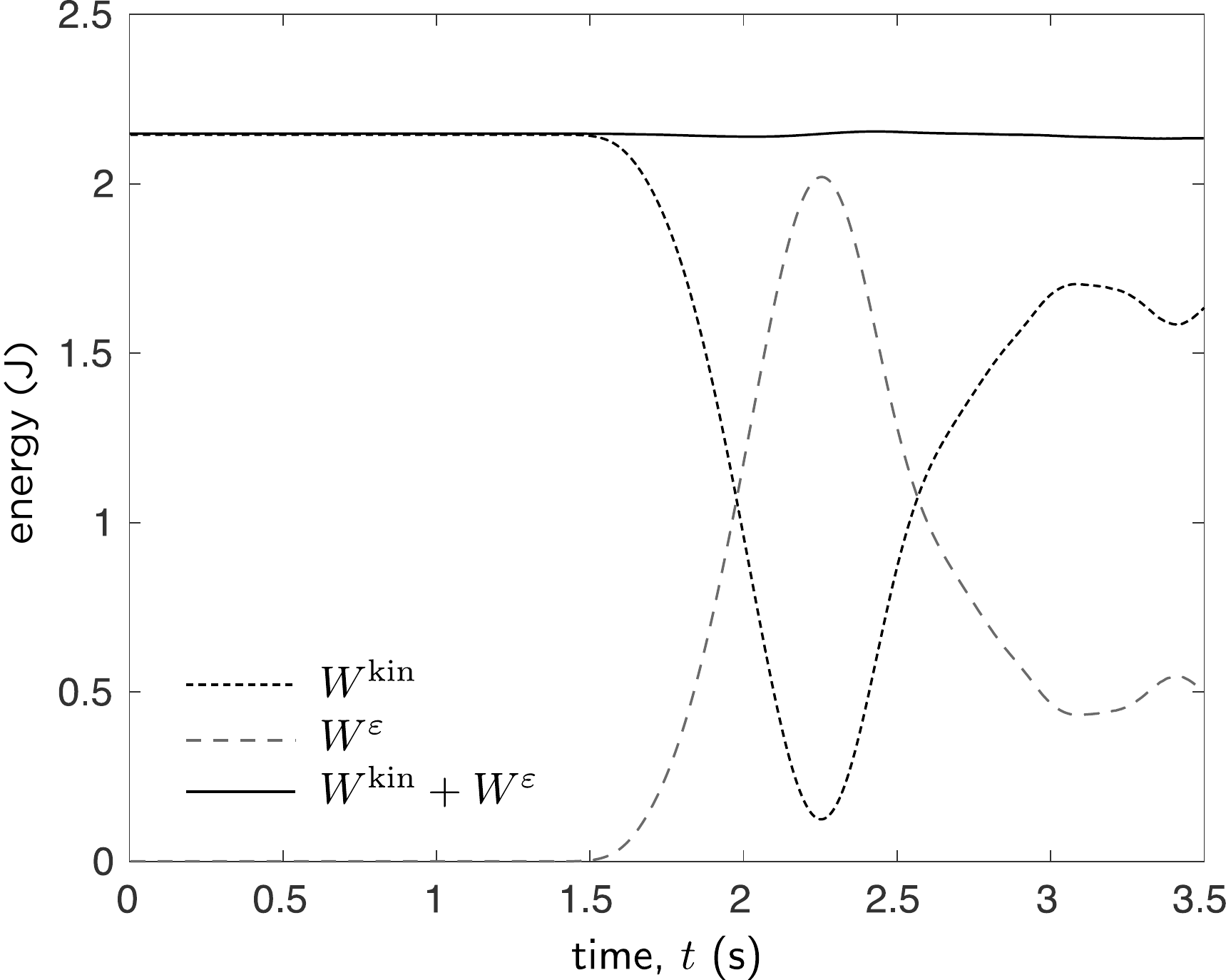}
        \caption{\footnotesize USF Ghost-stabilised energy variation }
    \end{subfigure} \hspace{0.01\textwidth}%
    \begin{subfigure}[c]{0.49\textwidth}
        \centering
        \includegraphics[width=\textwidth]{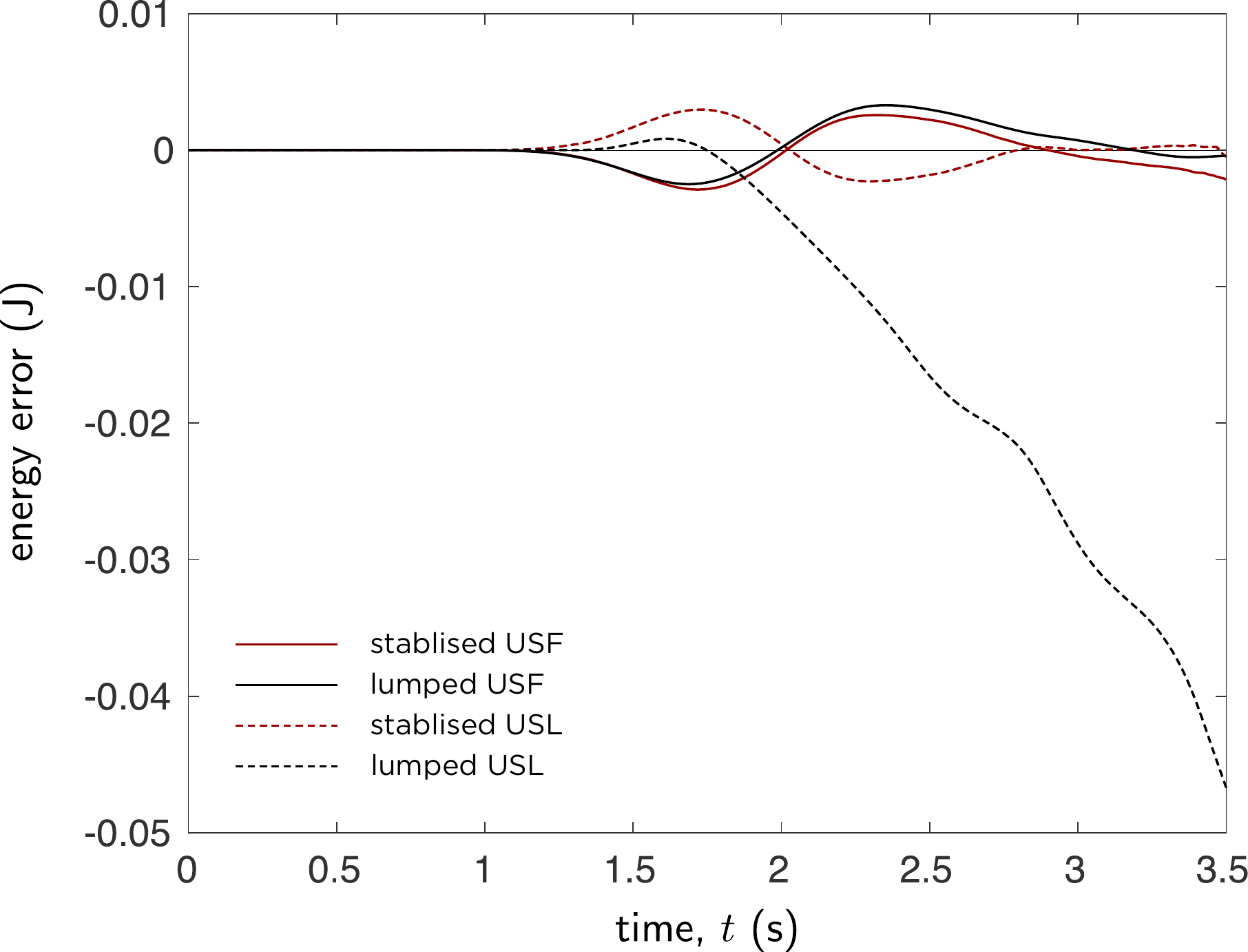}
        \caption{\footnotesize energy error}
    \end{subfigure}
\caption{Collision of elastic ghosts: (a) USF Ghost-stabilised energy variation  for $h=0.1$m, $n_t=1000$ and the \emph{medium} material point discretisation and (b) energy error for $h=0.1$m, $n_t=2000$ and the \emph{medium} material point discretisation.}
\label{Fig:ghostCollision}
\end{figure}

Figure~\ref{Fig:ghostCollision}~(a) shows the USF Ghost-stabilised energy variation  for $h=0.1$m, $n_t=1000$ and the \emph{medium} (2,389 material points for each ghost, see Section~\ref{sec:RBT}) material point discretisation and  Figure~\ref{Fig:ghostCollision}~(b) shows the energy error
\[
	\text{energy error} = W_i^{\text{kin}}+W_i^{\varepsilon}-W_0,
\]
where $W_i^{\text{kin}}$  and $W_i^{\varepsilon}$ are the kinetic and strain energies associated with time step $i$, for $h=0.1$m, $n_t=4000$ and the \emph{medium} material point discretisation over the duration of the analysis for the lumped and stabilised consistent mass matrices using an USF and USL approach.  As reported in the literature \cite{Bardenhagen2002,Love2006}, combining an USL approach with a lumped mass matrix (dashed black line) results in significant energy dissipation.  This excessive dissipation is removed via the adoption of the Ghost stabilised consistent mass matrix (dashed red line), resulting a simulation that approximately mirrors the USF analysis (solid red line) in terms of the energy dissipated at any point in the analysis.  The lumped and stabilised consistent mass matrices yield very similar results in terms of the energy error for the USF approach.

Table~\ref{Tab:ElasticGhostErrors} provides the normalised mean energy error
\[
	\text{normalised mean error} = \frac{1}{n_{t} W_0} \sum_{i=1}^{n_t} \bigl|W_i^{\text{kin}}+W_i^{\varepsilon}-W_0 \bigr|
\]
for the 108 analyses, where $n_t$ is the total number of time steps.  The numbers are coloured according to the magnitude of the average energy error, when yellow and green correspond to the minimum and maximum errors, respectively.  A number of other points can be observed:
\begin{enumerate}[(i)]  
	\item increasing the number of time steps reduces the average energy error over the simulation for all analyses;
	\item the Ghost stabilised consistent mass matrix is more sensitive to changes in background grid resolution compared to the lumped mass matrix when an USF approach is adopted, with the average energy error being smaller for the $h=0.2$m and $0.1$m analyses when using $[\tilde{M}_v]$ than $[\bar{M}_v]$ but the opposite being true for the $h=0.05$m analysis;
	\item the combination of an USL approach with a lumped mass matrix, $[\bar{M}_v]$, yields excessive energy dissipation, as shown by sub Tables~\ref{Tab:ElasticGhostErrors}~(j), (k) and (l), where the average energy error an order of magnitude larger than that of the USL with a stabilised consistent mass matrix, $[\tilde{M}_v]$, for mesh sizes $h=0.1$ and $0.05$m; 
\end{enumerate} 
The analyses presented for this problem demonstrate that the Ghost stabilisation allows the consistent mass matrix to be used for explicit dynamic analysis and also corrects an issue with the USL approach when combined with a lumped mass matrix.  The following section will analyse the impact of including elasto-plasticity.  

\begin{table}[h!]
	\footnotesize
        \begin{center}
	\caption{Collision of elastic ghosts: normalised mean energy errors ($\times10^{-3}$) for different discretisations using USL and USF approaches and the lumped and stabilised consistent mass matrices.  The cells are coloured according to the magnitude of the error from yellow (small) to green (large).}
	\label{Tab:ElasticGhostErrors}
	\begin{subtable}[h]{0.32\textwidth}
		\centering 
		\caption{\scriptsize USF ghost coarse}
            \begin{tabular}{l l | R | R | R}
	 \multicolumn{2}{c|}{$h$ (m)} & \multicolumn{1}{c|}{$0.20$} & \multicolumn{1}{c|}{$0.10$} & \multicolumn{1}{c}{$0.05$} \\  \cline{1-4} 
	\multirow{3}{*}{\rotatebox[origin=c]{90}{time steps~}~} 
	&1000~~& 1.370 & 0.888 &  1.386\\
	&2000~~& 0.682 & 0.445 &   0.707\\
	&4000~~& 0.340 & 0.223 &   0.357\\ 
            \end{tabular}\vspace*{2mm}
	\end{subtable} 
	\begin{subtable}[h]{0.32\textwidth}
		\centering 
		\caption{\scriptsize USF ghost medium}
            \begin{tabular}{l l | R | R | R}
	 \multicolumn{2}{c|}{$h$ (m)} & \multicolumn{1}{c|}{$0.20$} & \multicolumn{1}{c|}{$0.10$} & \multicolumn{1}{c}{$0.05$} \\  \cline{1-4} 
	\multirow{3}{*}{\rotatebox[origin=c]{90}{time steps~}~} 
	&1000~~& 1.294 &  0.889  &  1.559\\
	&2000~~& 0.643 &   0.447  &  0.792\\
	&4000~~& 0.321 &  0.224  &  0.399\\ 
            \end{tabular}\vspace*{2mm}
	\end{subtable} 
	\begin{subtable}[h]{0.32\textwidth}
		\centering 
		\caption{\scriptsize USF ghost fine}
            \begin{tabular}{l l | R | R | R}
	 \multicolumn{2}{c|}{$h$ (m)} & \multicolumn{1}{c|}{$0.20$} & \multicolumn{1}{c|}{$0.10$} & \multicolumn{1}{c}{$0.05$} \\  \cline{1-4} 
	\multirow{3}{*}{\rotatebox[origin=c]{90}{time steps~}~} 
	&1000~~& 1.261 &  0.899  &  1.629\\
	&2000~~& 0.627 &   0.451  &  0.828\\
	&4000~~& 0.313 &  0.226  &  0.418\\ 
            \end{tabular}\vspace*{2mm}
	\end{subtable}\\

	\begin{subtable}[h]{0.32\textwidth}
		\centering 
		\caption{\scriptsize USF lumped coarse}
            \begin{tabular}{l l | R | R | R}
	 \multicolumn{2}{c|}{$h$ (m)} & \multicolumn{1}{c|}{$0.20$} & \multicolumn{1}{c|}{$0.10$} & \multicolumn{1}{c}{$0.05$} \\  \cline{1-4} 
	\multirow{3}{*}{\rotatebox[origin=c]{90}{time steps~}~} 
	&1000~~& 1.418 &   0.956 &   0.871\\
	&2000~~& 0.706 &   0.478  &  0.438\\
	&4000~~& 0.352 &  0.239  &  0.220\\ 
            \end{tabular}\vspace*{2mm}
	\end{subtable} 
	\begin{subtable}[h]{0.32\textwidth}
		\centering 
		\caption{\scriptsize USF lumped medium}
            \begin{tabular}{l l | R | R | R}
	 \multicolumn{2}{c|}{$h$ (m)} & \multicolumn{1}{c|}{$0.20$} & \multicolumn{1}{c|}{$0.10$} & \multicolumn{1}{c}{$0.05$} \\  \cline{1-4} 
	\multirow{3}{*}{\rotatebox[origin=c]{90}{time steps~}~} 
	&1000~~& 1.325 &  0.916  & 1.011\\
	&2000~~& 0.660 &  0.458  &  0.511\\
	&4000~~& 0.329 &  0.229  & 0.257\\ 
            \end{tabular}\vspace*{2mm}
	\end{subtable} 
	\begin{subtable}[h]{0.32\textwidth}
		\centering 
		\caption{\scriptsize USF lumped fine}
            \begin{tabular}{l l | R | R | R}
	 \multicolumn{2}{c|}{$h$ (m)} & \multicolumn{1}{c|}{$0.20$} & \multicolumn{1}{c|}{$0.10$} & \multicolumn{1}{c}{$0.05$} \\  \cline{1-4} 
	\multirow{3}{*}{\rotatebox[origin=c]{90}{time steps~}~} 
	&1000~~& 1.283 &   0.908 &   1.122\\
	&2000~~& 0.638 &    0.454  &  0.567\\
	&4000~~& 0.319 &   0.227  &  0.285\\ 
            \end{tabular}\vspace*{2mm}
	\end{subtable}\\

	\begin{subtable}[h]{0.32\textwidth}
		\centering 
		\caption{\scriptsize USL ghost coarse}
            \begin{tabular}{l l | R | R | R}
	 \multicolumn{2}{c|}{$h$ (m)} & \multicolumn{1}{c|}{$0.20$} & \multicolumn{1}{c|}{$0.10$} & \multicolumn{1}{c}{$0.05$} \\  \cline{1-4} 
	\multirow{3}{*}{\rotatebox[origin=c]{90}{time steps~}~} 
	&1000~~& 4.533 &   0.859  &  0.980\\
	&2000~~& 2.738 &   0.493  &  0.557\\
	&4000~~& 1.542 & 0.272  &  0.409\\ 
            \end{tabular}\vspace*{2mm}
	\end{subtable} 
	\begin{subtable}[h]{0.32\textwidth}
		\centering 
		\caption{\scriptsize USL ghost medium}
            \begin{tabular}{l l | R | R | R}
	 \multicolumn{2}{c|}{$h$ (m)} & \multicolumn{1}{c|}{$0.20$} & \multicolumn{1}{c|}{$0.10$} & \multicolumn{1}{c}{$0.05$} \\  \cline{1-4} 
	\multirow{3}{*}{\rotatebox[origin=c]{90}{time steps~}~} 
	&1000~~& 4.319 &  0.751 &   1.242\\
	&2000~~& 2.576 &   0.346  &  0.519\\
	&4000~~& 1.435 &  0.174 &   0.216\\ 
            \end{tabular}\vspace*{2mm}
	\end{subtable}
	\begin{subtable}[h]{0.32\textwidth}
		\centering 
		\caption{\scriptsize USL ghost fine}
            \begin{tabular}{l l | R | R | R}
	 \multicolumn{2}{c|}{$h$ (m)} & \multicolumn{1}{c|}{$0.20$} & \multicolumn{1}{c|}{$0.10$} & \multicolumn{1}{c}{$0.05$} \\  \cline{1-4} 
	\multirow{3}{*}{\rotatebox[origin=c]{90}{time steps~}~} 
	&1000~~& 4.183 &   0.789   &  1.425\\
	&2000~~& 2.440 &   0.368   &   0.639\\
	&4000~~& 1.337 &   0.175   &   0.275\\ 
            \end{tabular}\vspace*{2mm}
	\end{subtable}\\

	\begin{subtable}[h]{0.32\textwidth}
		\centering 
		\caption{\scriptsize USL lumped coarse}
            \begin{tabular}{l l | R | R | R}
	 \multicolumn{2}{c|}{$h$ (m)} & \multicolumn{1}{c|}{$0.20$} & \multicolumn{1}{c|}{$0.10$} & \multicolumn{1}{c}{$0.05$} \\  \cline{1-4} 
	\multirow{3}{*}{\rotatebox[origin=c]{90}{time steps~}~} 
	&1000~~& 12.32 &   8.417  & 10.58\\
	&2000~~& 6.647 &   4.726  &   6.790\\
	&4000~~& 3.460 &  2.521  &  3.989\\ 
            \end{tabular}\vspace*{2mm}
	\end{subtable} 
	\begin{subtable}[h]{0.32\textwidth}
		\centering 
		\caption{\scriptsize USL lumped medium}
            \begin{tabular}{l l | R | R | R}
	 \multicolumn{2}{c|}{$h$ (m)} & \multicolumn{1}{c|}{$0.20$} & \multicolumn{1}{c|}{$0.10$} & \multicolumn{1}{c}{$0.05$} \\  \cline{1-4} 
	\multirow{3}{*}{\rotatebox[origin=c]{90}{time steps~}~} 
	&1000~~& 12.39 &   8.207 &   8.899\\
	&2000~~&  6.684 &   4.607 &   5.674\\
	&4000~~&  3.481 &  2.457  &  3.319\\ 
            \end{tabular}\vspace*{2mm}
	\end{subtable}
	\begin{subtable}[h]{0.32\textwidth}
		\centering 
		\caption{\scriptsize USL lumped fine}
            \begin{tabular}{l l | R | R | R}
	 \multicolumn{2}{c|}{$h$ (m)} & \multicolumn{1}{c|}{$0.20$} & \multicolumn{1}{c|}{$0.10$} & \multicolumn{1}{c}{$0.05$} \\  \cline{1-4} 
	\multirow{3}{*}{\rotatebox[origin=c]{90}{time steps~}~} 
	&1000~~& 12.35 &    7.992  &  8.022\\
	&2000~~& 6.666 &    4.484   & 5.105\\
	&4000~~& 3.473 &   2.391 &   2.982 \\ 
            \end{tabular}\vspace*{2mm}
	\end{subtable}
        \end{center}
\end{table}

\subsection{Explicit dynamics: collision of elasto-plastic ghosts}

This example considers the collision of two elasto-plastic ghosts to demonstrate that the Ghost stabilisation approach can be adopted for non-linear materials.  The analyses were run on a 1m$\times$1m background grid, with 20 elements in each direction ($h=0.05$m). The ghosts were initially placed in the upper right and lower left corners of the background grid, with velocities of $\{-0.1,-0.1\}$m/s and $\{0.1, 0.1\}$m/s respectively (as shown in Figure~\ref{Fig:ghostCollisionSetup}). Each ghost had width and height of $0.4$m, a Young's modulus of $E=1000$Pa, a Poisson's ratio of $\nu=0.3$ and a density of $\rho=1000$kg/m\textsuperscript{3}.  Yielding of the material was governed by a von Mises function of the form
\begin{equation}\label{eqn:vM}
  f = \rho-\rho_y=0,
\end{equation} 
where $\rho=\sqrt{2J_2}$, $J_2=s_{ij}s_{ji}/2$, $s_{ij}=\tau_{ij}-\tau_{kk}\delta_{ij}/3$ and $\tau_{ij}$ is the Kirchhoff stress.  The deviatoric yield stress of both Ghosts was set to $\rho_y=100$Pa and the constitutive model adopted a elastic predictor, plastic corrector closest point projection algorithm (see \cite{coombs2011finitethesis} amongst others).  The \emph{medium} generalised interpolation material point resolution used in the previous sections, and detailed in Table~\ref{Tab:Ghost_noMPs}, was adopted for this analysis.  The effect of gravity was ignored and the behaviour modelled over 3.5 seconds using 1000 time steps. The system was undamped and was analysed using the USF explicit dynamic material point formulation described in Section~\ref{sec:ExpDyn} with a stabilised consistent mass matrix ($\gamma_M=\frac{1}{4}\rho$).

\begin{figure}[!h]
\centering
\includegraphics[width=0.6\textwidth]{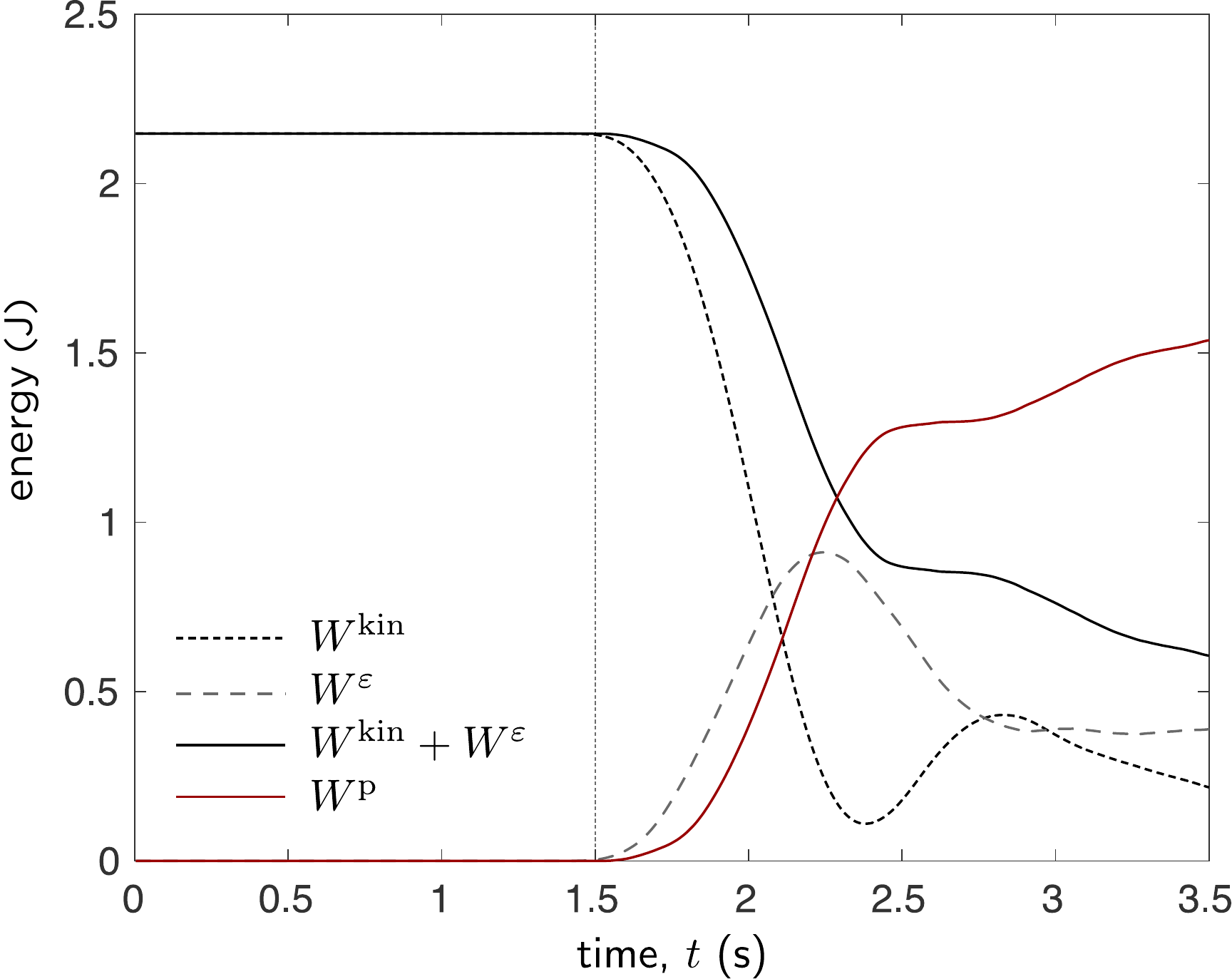}
\caption{Collision of elasto-plastic ghosts: energy components for a \emph{medium} material point discretisation, a background mesh size of $h=0.05$m and 1000 time steps.}
\label{Fig:EPghostEnergy}
\end{figure}

Figure~\ref{Fig:EPghostEnergy} shows the evolution of the energy components over the analysis, showing the  kinetic ($W^{\text{kin}}$, dashed black line) and strain ($W^{\varepsilon}$, dashed grey line) energies, the plastic dissipation ($W^{\text{p}}$, red line) and the total energy of the bodies ($W^{\text{kin}}+W^{\varepsilon}$, black line).  The initial collision starts at approximately $1.5$s, which initiates strain energy generation and plastic dissipation within the ghosts.  Beyond this initial collosion there is additional plastic dissipation from around 2.7s due to stress waves through the material.  The plastic deformation also \emph{traps} a degree of strain energy in the ghosts, as seen by the plateau in the strain energy beyond 3.0s.

\begin{figure}[!h]
\centering
    \begin{subfigure}[t]{0.35\textwidth}
        \centering
        \includegraphics[width=\textwidth]{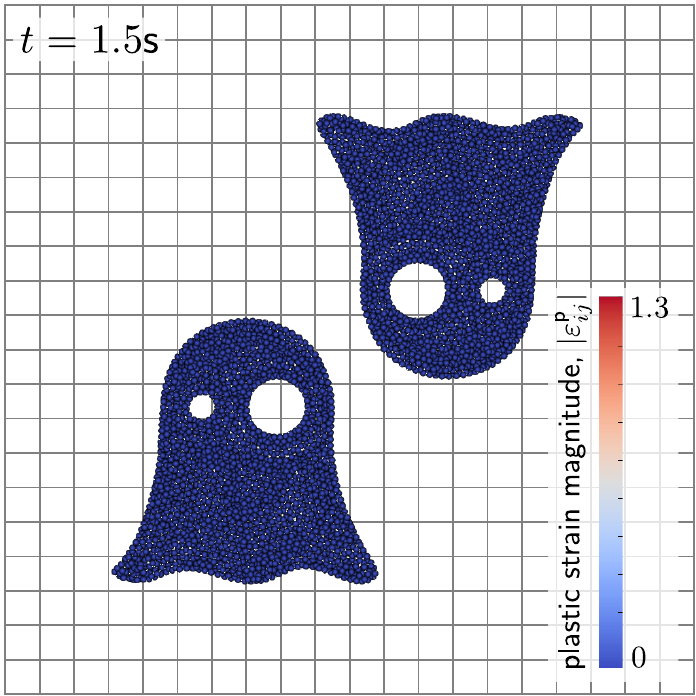}
        \caption{\footnotesize $t=1.5$s}
    \end{subfigure}\hspace*{5mm}%
    \begin{subfigure}[t]{0.35\textwidth}
        \centering
        \includegraphics[width=\textwidth]{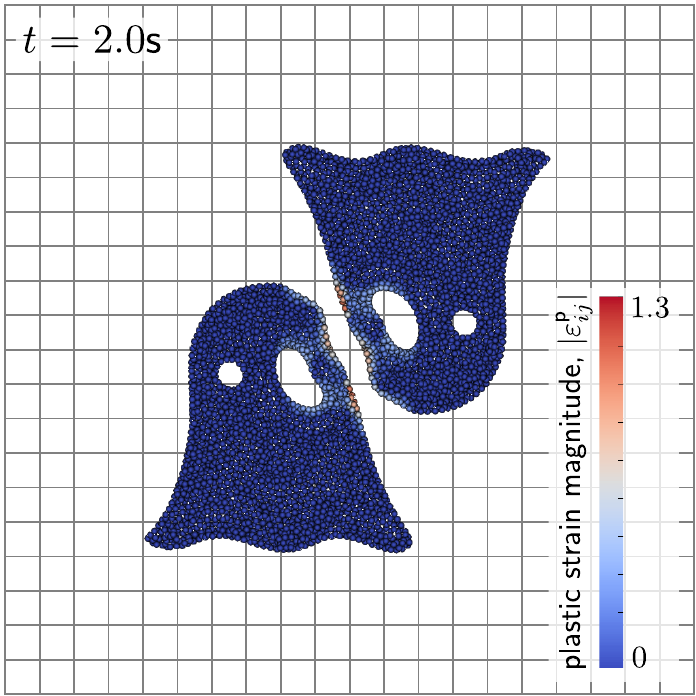}
        \caption{\footnotesize $t=2.0$s}
    \end{subfigure}\\[2mm]
    \begin{subfigure}[t]{0.35\textwidth}
        \centering
        \includegraphics[width=\textwidth]{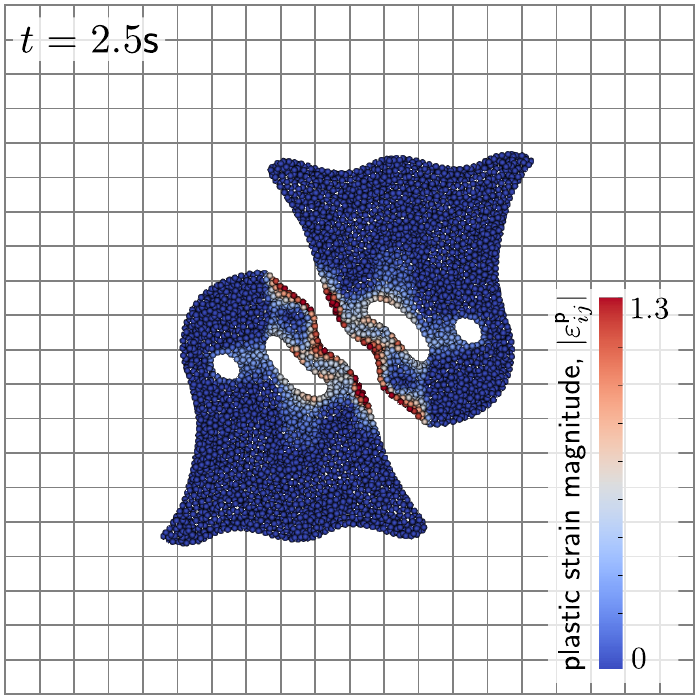}
        \caption{\footnotesize $t=2.5$s}
    \end{subfigure}\hspace*{5mm}%
     \begin{subfigure}[t]{0.35\textwidth}
        \centering
        \includegraphics[width=\textwidth]{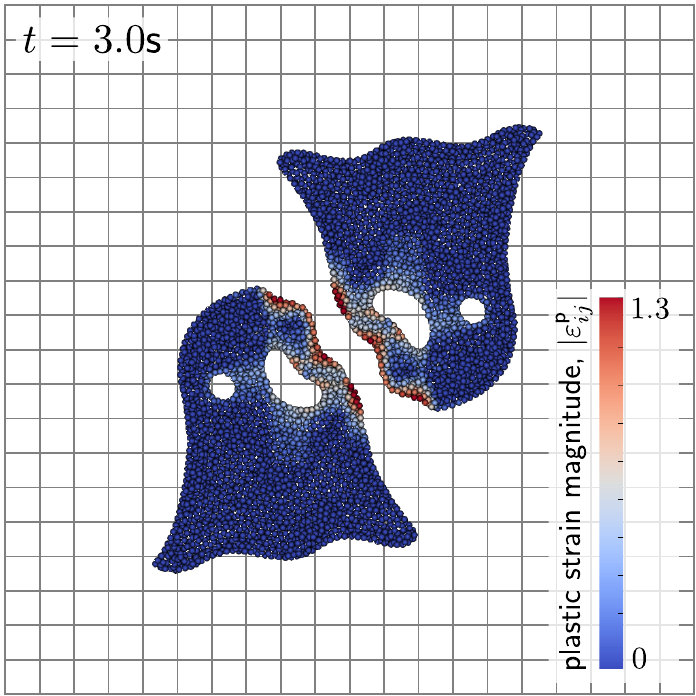}
        \caption{\footnotesize $t=3.0$s}
    \end{subfigure}
\caption{Collision of elasto-plastic ghosts: deformed material point positions coloured according to the magnitude of the total plastic logarithmic strains using a \emph{medium} material point discretisation, a background mesh size of $h=0.05$m and 1000 time steps.}
\label{Fig:EPghostsPlasticStrains}
\end{figure}

Figure~\ref{Fig:EPghostsPlasticStrains} shows the material point positions at $t=1.5$s, $2.0$s, $2.5$s and $3.0$s, coloured according to the magnitude of the total plastic logarithmic strains.  It is worth noting that the native contact in material point methods is based on the background grid and the proximity between the two bodies when contact occurs is based on the grid spacing rather than a description of the surface of the objects.  This is why the bodies always maintain some separation.  Several papers (see \citet{Acosta2021} for a recent overview) have investigated alternative contact formulations in the material point method but these techniques are not explored in this paper, where the focus is on the stabilisation of material point methods in general, including both implicit \emph{quasi}-static analysis which will be explored in the next example.

\subsection{Implicit \emph{quasi}-static: elastic compression under self weight}

\begin{figure}[!h]
        \centering
        \includegraphics[width=0.8\textwidth]{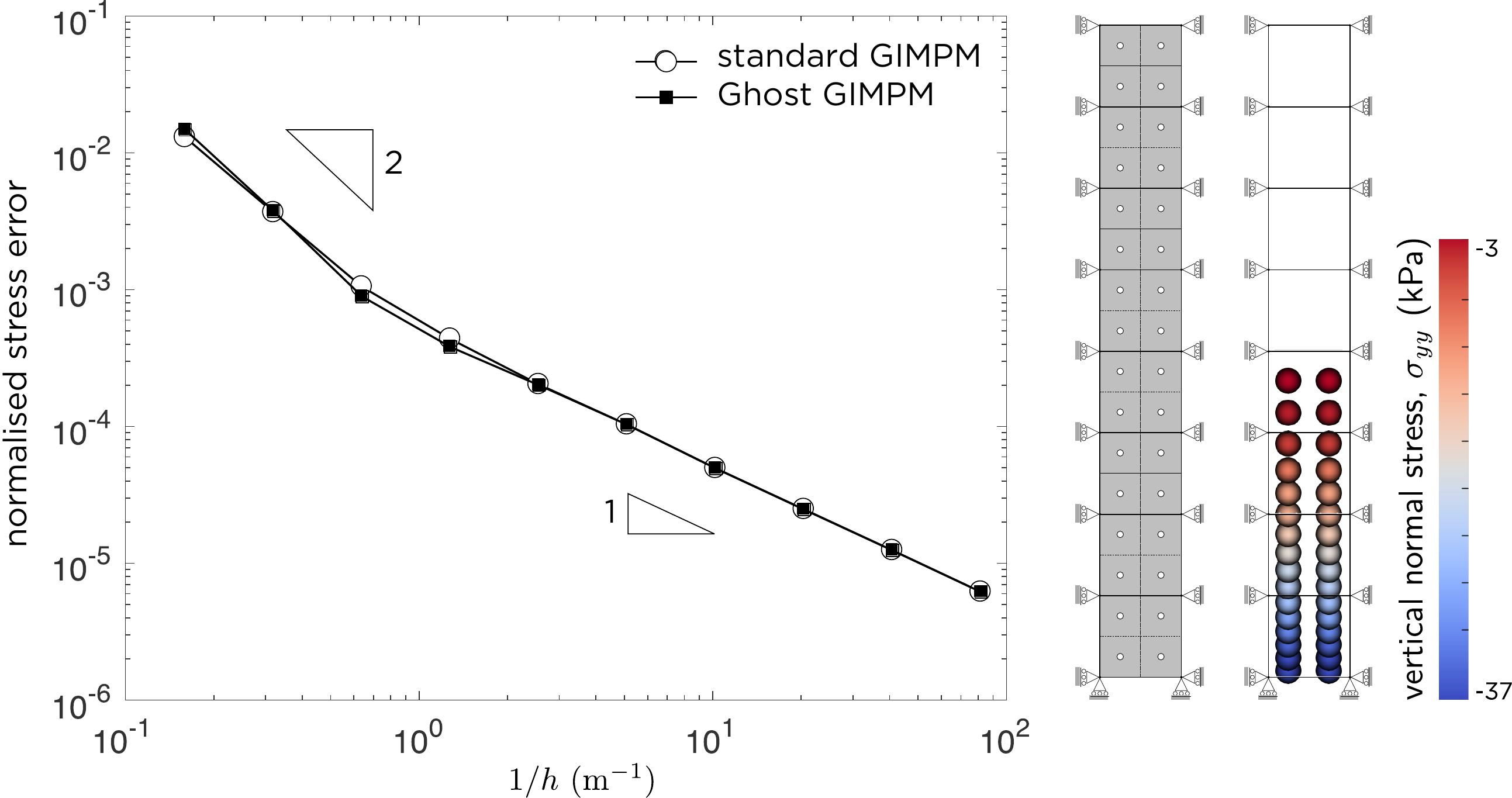}
	\caption{Compression under self weight: normalised stress error convergence with grid refinement and initial problem discretisation \& final material point positions coloured according to $\sigma_{yy}$ with $h=6.25$m.}
	\label{Fig:compactConverg}
\end{figure}

This example considers the one dimensional compression of an elastic column with an initial height of $l_0=50$m under its own self weight and aims to demonstrate that including Ghost stabilisation has minimal impact on the convergence of the method using a problem with an analytical solution.  The material has a Young's modulus of 10kPa and a Poisson's ratio of zero.  The background mesh is comprised of square  background elements with roller boundary conditions on the base and sides and the column is discretised by a $2\times2$ grid of equally spaced material points in each initially populated background grid element (as shown to the right of Figure~\ref{Fig:compactConverg} for $h=6.25$m).  A body force of 800N/m$^2$ ($g=10$m/s$^2$ and an initial density of $\varrho_0=80$kg/m$^3$) is applied over 40 equal load steps.  The magnitude of the load causes the column to compress to approximately half of its initial height.   

The analytical solution for the normal stress in the vertical ($y$) direction for this problem is
\begin{equation}
  \sigma^a_{yy} =  \varrho_0 g (l_0-Y),
\end{equation}
where $Y$ is the original position of the point in the body and $l_0$ is the original height of the column.  Figure~\ref{Fig:compactConverg} shows the convergence of the generalised interpolation material point method with background mesh refinement whilst maintaining $4$ material points per initially populated element.  The reported normalised stress error is 
\begin{equation}
  \text{error} = \sum_{p=1}^{n_{p}}\frac{||\sigma^{p}_{yy}-\sigma^{a}_{yy}(Y_{p})||~V^0_p}{(g \varrho_0 l_0) V_0},
\end{equation}
where $V_0=\sum V^0_p$ is the initial volume of the column and $\sigma^{p}_{yy}$ is the vertical stress at each of the material points.  Figure~\ref{Fig:compactConverg} provides the errors with and without Ghost stabilisation, where $\gamma_k$ was set to $10$kPa for the stabilised analyses.  The generalised interpolation material point method with and without stabilisation converge towards the analytical solution at a rate between 1 and 2, which is consistent with the underlying basis of the method and there is very little difference between the error values for different background mesh resolutions.  The difference in error also reduces with grid refinement due to the region influenced by the stabilisation reducing as the grid is refined.  This point is reinforced by Table~\ref{Tab:compact}, which gives the normalised stress error for $h=50/2^9\approx 0.098$m and   $h=50/2^5\approx 0.163$m with different stabilisation parameter values.  It is clear from the table that reducing the background element size reduces the sensitivity of the result to the value of the penalty parameter for the reason explained above.  For the analysis with a mesh size of $0.098$m, twelve orders of magnitude variation in the penalty parameter only changes the normalised stress error by $0.49$\%, and the variation in stress error is 15\% for the $h=1.563$m analyses.  Although a variation of 15\% in the normalised stress error seems significant, it is only a very small variation in the actual error, which is magnified by normalising by the error associated with the non-stabilised analysis.


\begin{table}[!h]
\centering
\caption{Compression under self weight: stress error variation with $\gamma_k$.}
{\small \begin{tabular}{c | ccccc}
	& \multicolumn{5}{c}{normalised penalty parameter, $\gamma_k/E$}\\
	$h$ (m)~~  & ~~$1^{-6}$~~ & ~~$1^{-3}$~~ & ~~$1^0$~~ & ~~$1^{3}$~~ & ~~$1^{6}$~~ \\ \hline 
	$0.098$~~ & ~~$4.938\times10^{-5}$~~ & ~~$4.938\times10^{-5}$~~ & ~~$4.938\times10^{-5}$~~ & ~~$4.926\times10^{-5}$~~ & ~~$4.962\times10^{-5}$~~ \\
	$1.563$~~ & ~~$1.049\times10^{-3}$~~ & ~~$1.042\times10^{-3}$~~ & ~~$8.912\times10^{-4}$~~ & ~~$1.208\times10^{-3}$~~ & ~~$1.210\times10^{-3}$~~ 
\end{tabular}}
\label{Tab:compact}
\end{table}

\subsection{Implicit \emph{quasi}-static: elastic beam}

\begin{figure}[!h]
\centering
\includegraphics[width=0.75\textwidth]{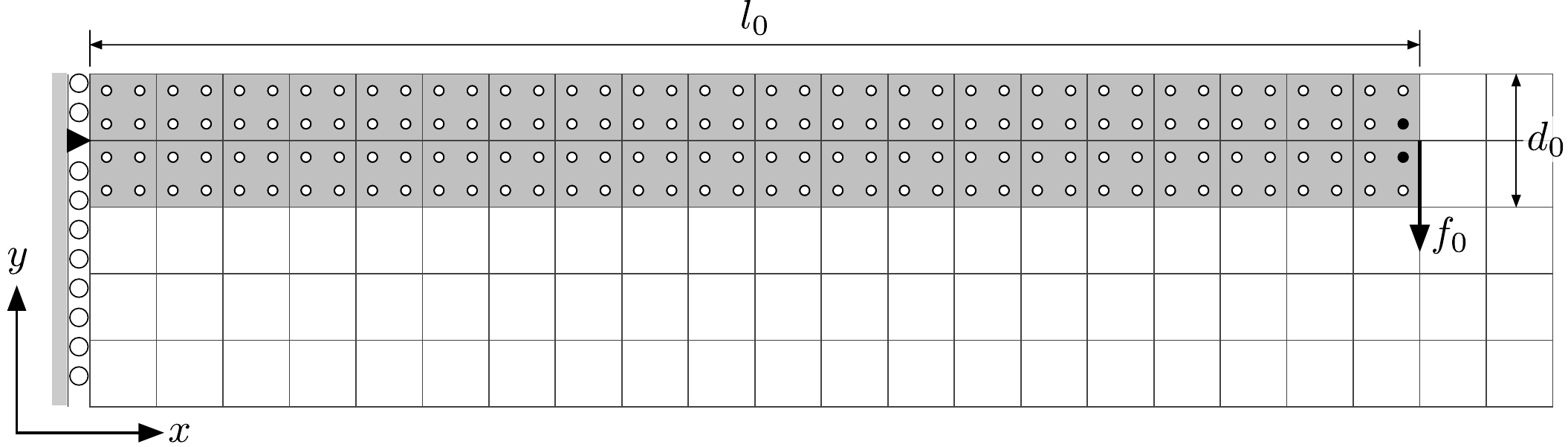}
\caption{elastic beam: initial discretisation and boundary conditions with $h=0.5$m (only part of the background mesh is shown).}
\label{Fig:beamSetup}
\end{figure}

The final example in this paper presents the implicit \emph{quasi}-static analysis of a large deformation beam as a challenging problem for the material point method due numerous small overlaps between the material points and the background grid.  The elastic cantilever beam was subjected to a point load at its free end and modelled using the generalised interpolation material point method.  The beam was $l_0=10$m long and $d_0=1$m deep and the material had a Young's modulus of $12$MPa, a Poisson's ratio of $0.2$ and was assumed to be weightless ($\rho=0$).    The $f_0=100$kN end point load was split between the two material points closest to the end of the beam either side of the neutral axis and applied over 50 equal load steps.  The initial discretisation of the beam is shown in Figure~\ref{Fig:beamSetup} with $h=0.5$m and with $2^2$ material points per initially populated background grid cells.  The loaded material points are shown by the black-filled circles.   The stiffness stabilisation parameter was taken to be the same as the Young's modulus of the material, $\gamma_K=12$MPa.

Table~\ref{Tab:beam} provides information on the stability of the \emph{quasi}-static analysis with different numbers of material point per element with and without Ghost stabilisation.  The table reports different information depending on if an analysis was able to complete all load steps:
\begin{itemize}
	\item[\xmark] unstable analysis: the number indicates the final stable load step; and
	\item[\cmark] stable analysis: the numbers indicate the total and (maximum) number of Newton-Raphson iterations for the analysis.
\end{itemize}  
The normalised convergence tolerance\footnote{The convergence criteria used in the implicit \emph{quasi}-static material point method implementation adopted in this paper is the same as that used by \citet{Coombs2020}, where the residual out of balance force (the difference between the internal force and the external actions), normalised by the magnitude of the external actions is checked until it converges below a given tolerance.  The only change for this paper is that the residual force equation now includes a contribution from the Ghost stabilisation, as explained in Section~\ref{sec:stiffStab}.} was set to $1\times10^{-6}$ for all analyses.  All of the Ghost-stabilised analyses were able to complete all of the load steps, however several of the standard GIMPM analyses failed to converge (or reached the maximum number of Newton iterations, set to $10$ for all analyses) before the full load had been applied.  For the standard GIMPM there is a general trend of increasing stability with increasing numbers of material points but this is not always the case, for example the $h=0.5$m, $6^2$ material point failed at load step 10 whereas the $5^2$ material point analysis was able to apply the full load for the same background mesh size.  In essence a user does not know if the analysis will converge or not before running a simulation when using the standard GIMPM, even if \emph{similar} analyses have converged.  Reducing the mesh size also increases the average total number of Newton-Raphson iterations required to complete the analysis for the standard (non-stabilised) GIMPM whereas the total and maximum number of iterations for the stabilised GIMPM is relatively insensitive to mesh size variations.  In fact, there is trend of reducing total numbers of iterations as the mesh is refined.         

\begin{table}[!h]
\centering
\caption{Elastic beam: for entries marked \xmark, the number indicates the final stable load step; for entries marked \cmark, the numbers indicate the total and (maximum) number of Newton iterations for the analysis.}
{\small \begin{tabular}{c | c c c | c c c}
	& \multicolumn{3}{c|}{standard GIMPM} & \multicolumn{3}{c}{Ghost stabilised GIMPM} \\
	MPs/elem. & ~~$h=0.500$~~ & ~~$h=0.250$~~ & ~~$h=0.125$~~ & ~~$h=0.500$~~ & ~~$h=0.250$~~ & ~~$h=0.125$~~ \\ \hline 
	$2^2$ 	& \xmark~5 		& \xmark~11 		& \xmark~1	& \cmark~205(5) 	& \cmark~203(5) & \cmark~205(5)\\
	$3^2$	& \cmark~215(9) 	& \cmark~215(8)	& \xmark~5	& \cmark~204(5)	& \cmark~204(5) & \cmark~203(5)\\
	$4^2$	& \xmark~40 		& \xmark~26		& \cmark~215(5)& \cmark~204(5)	& \cmark~204(5) & \cmark~201(5)\\
	$5^2$	& \cmark~207(6) 	& \cmark~207(5)	& \cmark~214(5)& \cmark~204(5)	& \cmark~203(5) & \cmark~202(5)\\
	$6^2$	& \xmark~9		& \cmark~207(5)	& \cmark~216(6)& \cmark~204(5)	& \cmark~204(5) & \cmark~202(5)\\ \hline
	average 	& ~~211(7.5) 		&  ~~209.7(6) 		& ~~215(5.3)	& ~~204.2(5) 		& ~~203.6(5) 	 & ~~202.6(5) 
\end{tabular}}
\label{Tab:beam}
\end{table}

\begin{figure}[!h]
\centering
\includegraphics[width=0.7\textwidth]{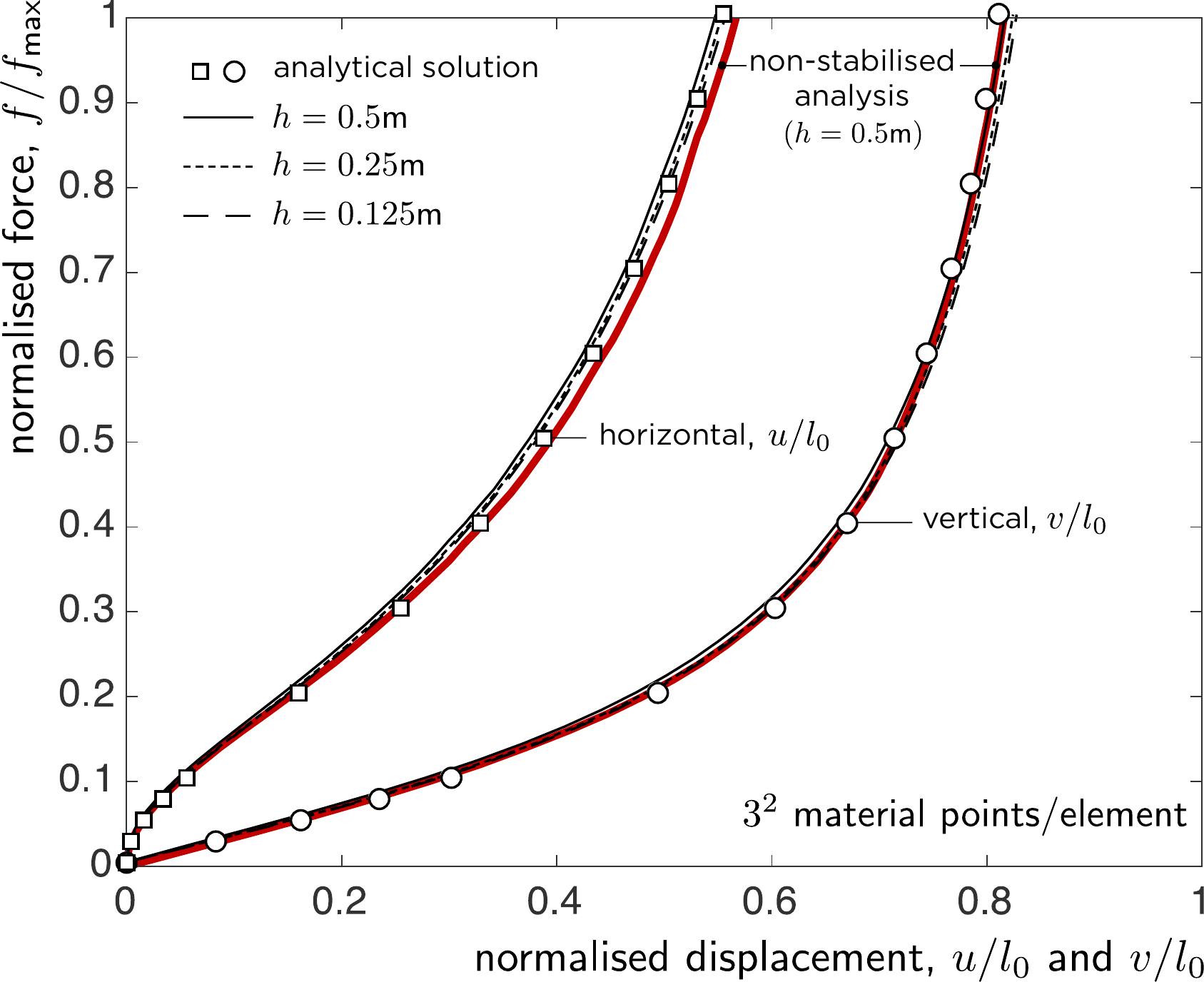}
\caption{elastic beam: normalised force versus displacement (where the squares and circles show the analytical solution of \citet{Molstad1977} for the horizontal and vertical displacements, respectively)  with $3^2$ material points per initially populated element.  The red line shows the non-stabilised result with $h=0.5$m. }
\label{Fig:beam}
\end{figure}

Figure~\ref{Fig:beam} shows the normalised force-displacement response of the beam for the Ghost stabilised GIMPM with different background mesh sizes and $3^2$ material points per initially populated element.  The squares and circles show the analytical solution of \citet{Molstad1977} for the horizontal and vertical displacements, respectively, and the Ghost GIMPM results are shown by the black lines.  All of the Ghost GIMPM results are in good agreement with the analytical solution.  The simulations do over estimate the vertical displacement of the beam compared to the analytical solution but this is consistent with other numerical results in the literature \cite{Coombs2020,Charlton2017,Coombs2020a} and is linked to the analytical assumption that the beam does not change in length.  The red line in Figure~\ref{Fig:beam} shows the non-stabilised GIMPM result with $h=0.5$m, which deviates from the expected (analytical) response for the horizontal displacement, showing under-stiff behaviour. 


\begin{figure}[!h]
\centering
\includegraphics[width=\textwidth]{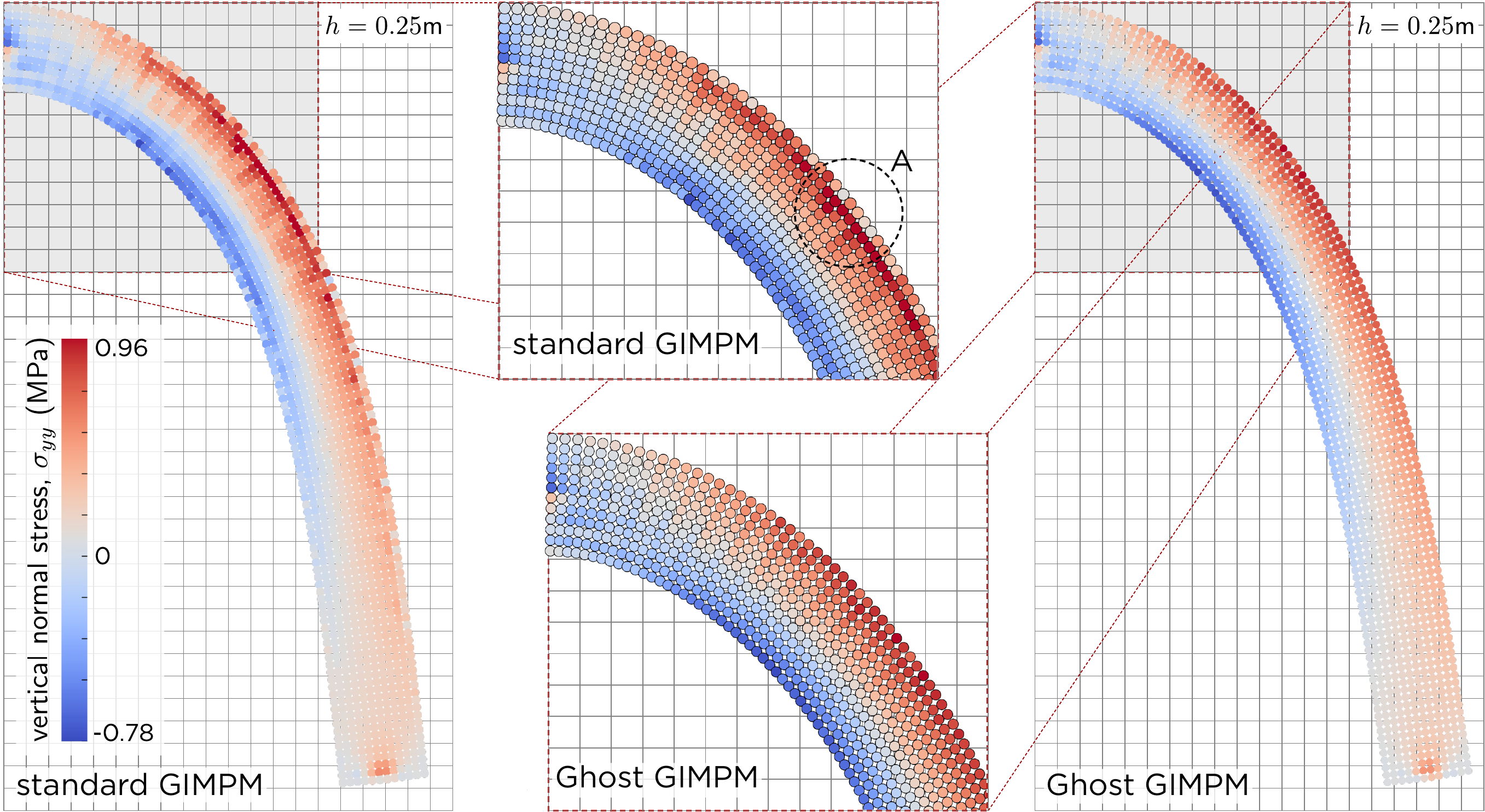}
\caption{elastic beam: vertical normal stress, $\sigma_{yy}$, distribution at the end of the analysis for the non-stabilised and Ghost stabilised GIMPMs with $h=0.25$m and $3^2$ material point per initially populated background grid cell. }
\label{Fig:beamStress}
\end{figure}

The vertical normal stress, $\sigma_{yy}$, distributions at the end of the analysis for the standard, non-stabilised (left) and Ghost stabilised (right) GIMPMs with $h=0.25$m and $3^2$ material point per initially populated background grid cell are shown in Figure~\ref{Fig:beamStress}.  The key difference in the stress distributions between the two methods is on the top and bottom surfaces of the beams, where the standard GIMPM predicts spurious oscillations in the stress field.  This variation in stress is particularly evident around region {\sf A} within the top inset figure, where a drop in the tensile stress on the top surface of the beam is due to the poor distribution of material points within the background grid at the edge of the physical domain.  These stress oscillations are removed by the Ghost stabilisation due to the penalisation of variation in the gradient of the solution field over the boundary element edges.

\subsection{Observations}

This section has presented five numerical investigations for dynamic and \emph{quasi}-static large deformation stress analysis problems.  The following key observations can be drawn from these numerical examples:
\begin{enumerate}[(i)]
	\item for explicit dynamics the Ghost stabilisation technique: 
	\begin{enumerate}[(a)]
	\item opens the door to the consistent mass matrix being used for practical simulations with arbitrary positioning of the physical domain (material points) on the background grid; 
	\item unlike the use of the lumped mass matrix, adopting the stabilised consistent mass matrix maintains the material point velocity field for rigid body motions and linear deformation fields; and
	\item corrects the excessive dissipation seen when combining USL approaches with a lumped mass matrix, due to consistent mapping of the nodal velocity field to the material points when determining changes in the deformation field.      
	\end{enumerate}
	\item for \emph{quasi}-static implicit analysis the Ghost stabilisation technique:
	\begin{enumerate}[(a)]
		\item maintains the ability of the underlying numerical algorithm to converge towards problems with analytical solutions as the stabilisation only influences a local region of the physical body near  the problem domain boundary; 
		\item removes the stability uncertainty when modelling large deformation problems - without stabilisation it is difficult to say if a analysis will converge or not even if similar analyses are stable; 
		\item permits the use of lower numbers of material points whilst maintaining stability of the overall method; and
		\item significantly reduces the stress oscillations seen near the boundary of the physical domain, which are due to small overlaps rather than cell crossing errors in non-stabilised methods.
	\end{enumerate}
\end{enumerate}

\section{Conclusion}

This paper has proposed a new stabilisation technique for explicit dynamic and implicit \emph{quasi}-static material point methods that resolves the vast majority of the stability issues encountered by the method, whilst not damaging their underlying properties.  The approach can be applied to all material point method variants, including implicit dynamics, without requiring an explicit description of the boundary of the problem.  The technique offered in this paper is one step on the road to the material point method becoming a usable tool for practical engineering analyses.

\section*{Acknowledgements}

The author would like to acknowledge the contributions of the Computational Mechanics Research Node in the Department of Engineering of Durham University.  The research presented in this article has benefited from discussions with, and feedback from, Charles Augarde, Robert Bird, Nathan Gavin, Ted O'Hare and Giuliano Pretti. 

This work was supported by the Engineering and Physical Sciences Research Council [grant numbers \href{http://gow.epsrc.ac.uk/NGBOViewGrant.aspx?GrantRef=EP/W000970/1}{EP/W000970/1}, \href{http://gow.epsrc.ac.uk/NGBOViewGrant.aspx?GrantRef=EP/R004900/1}{EP/R004900/1} and \href{http://gow.epsrc.ac.uk/NGBOViewGrant.aspx?GrantRef=EP/N006054/1}{EP/N006054/1}].    All data created during this research are openly available at \href{https://collections.durham.ac.uk/}{collections.durham.ac.uk/} (specific DOI to be confirmed if/when the paper is accepted).  For the purpose of open access, the author has applied a Creative Commons Attribution (CC BY) licence to any Author Accepted Manuscript version arising.

\bibliographystyle{model1b-num-names}
\bibliography{iGIMPrefs,mpmGhost}

\begin{thebibliography}{38}
\expandafter\ifx\csname natexlab\endcsname\relax\def\natexlab#1{#1}\fi
\providecommand{\url}[1]{\texttt{#1}}
\providecommand{\href}[2]{#2}
\providecommand{\path}[1]{#1}
\providecommand{\DOIprefix}{doi:}
\providecommand{\ArXivprefix}{arXiv:}
\providecommand{\URLprefix}{URL: }
\providecommand{\Pubmedprefix}{pmid:}
\providecommand{\doi}[1]{\href{http://dx.doi.org/#1}{\path{#1}}}
\providecommand{\Pubmed}[1]{\href{pmid:#1}{\path{#1}}}
\providecommand{\bibinfo}[2]{#2}
\ifx\xfnm\relax \def\xfnm[#1]{\unskip,\space#1}\fi
\bibitem[{Andersen and Andersen(2010)}]{Andersen2010}
\bibinfo{author}{S.~Andersen}, \bibinfo{author}{L.~Andersen},
  \bibinfo{title}{Analysis of spatial interpolation in the material-point
  method}, \bibinfo{journal}{Computers \& Structures} \bibinfo{volume}{88}
  (\bibinfo{year}{2010}) \bibinfo{pages}{506 -- 518}.
\bibitem[{Bardenhagen(2002)}]{Bardenhagen2002}
\bibinfo{author}{S.~Bardenhagen}, \bibinfo{title}{Energy conservation error in
  the material point method for solid mechanics}, \bibinfo{journal}{Journal of
  Computational Physics} \bibinfo{volume}{180} (\bibinfo{year}{2002})
  \bibinfo{pages}{383 -- 403}.
\bibitem[{Bardenhagen and Kober(2004)}]{Bardenhagen2004}
\bibinfo{author}{S.~Bardenhagen}, \bibinfo{author}{E.~Kober},
  \bibinfo{title}{The generalized interpolation material point method},
  \bibinfo{journal}{Computer Modeling in Engineering and Sciences}
  \bibinfo{volume}{5} (\bibinfo{year}{2004}) \bibinfo{pages}{477--496}.
\bibitem[{Berzins(2022)}]{Berzins2022}
\bibinfo{author}{M.~Berzins}, \bibinfo{title}{Energy conservation and accuracy
  of some {MPM} formulations}, \bibinfo{journal}{Computational Particle
  Mechanics}  (\bibinfo{year}{2022}) \bibinfo{pages}{1--13}.
\bibitem[{Bing et~al.(2019)Bing, Cortis, Charlton, Coombs and
  Augarde}]{Bing2019}
\bibinfo{author}{Y.~Bing}, \bibinfo{author}{M.~Cortis},
  \bibinfo{author}{T.~Charlton}, \bibinfo{author}{W.~Coombs},
  \bibinfo{author}{C.~Augarde}, \bibinfo{title}{B-spline based boundary
  conditions in the material point method}, \bibinfo{journal}{Computers \&
  Structures} \bibinfo{volume}{212} (\bibinfo{year}{2019}) \bibinfo{pages}{257
  -- 274}.
\bibitem[{Brackbill and Ruppel(1986)}]{Brackbill1986}
\bibinfo{author}{J.~Brackbill}, \bibinfo{author}{H.~Ruppel},
  \bibinfo{title}{Flip: A method for adaptively zoned, particle-in-cell
  calculations of fluid flows in two dimensions}, \bibinfo{journal}{Journal of
  Computational Physics} \bibinfo{volume}{65} (\bibinfo{year}{1986})
  \bibinfo{pages}{314 -- 343}.
\bibitem[{Burgess et~al.(1992)Burgess, Sulsky and Brackbill}]{Burgess1992}
\bibinfo{author}{D.~Burgess}, \bibinfo{author}{D.~Sulsky},
  \bibinfo{author}{J.~Brackbill}, \bibinfo{title}{Mass matrix formulation of
  the {FLIP} particle-in-cell method}, \bibinfo{journal}{Journal of
  Computational Physics} \bibinfo{volume}{103} (\bibinfo{year}{1992})
  \bibinfo{pages}{1--15}.
\bibitem[{Burman(2010)}]{Burman2010}
\bibinfo{author}{E.~Burman}, \bibinfo{title}{Ghost penalty},
  \bibinfo{journal}{Comptes Rendus Mathematique} \bibinfo{volume}{348}
  (\bibinfo{year}{2010}) \bibinfo{pages}{1217--1220}.
\bibitem[{Burman et~al.(2018)Burman, Elfverson, Hansbo, Larson and
  Larsson}]{Burman2018}
\bibinfo{author}{E.~Burman}, \bibinfo{author}{D.~Elfverson},
  \bibinfo{author}{P.~Hansbo}, \bibinfo{author}{M.~Larson},
  \bibinfo{author}{K.~Larsson}, \bibinfo{title}{Shape optimization using the
  cut finite element method}, \bibinfo{journal}{Computer Methods in Applied
  Mechanics and Engineering} \bibinfo{volume}{328} (\bibinfo{year}{2018})
  \bibinfo{pages}{242--261}.
\bibitem[{Chandra et~al.(2021)Chandra, Singer, Teschemacher, W\"{u}cheer and
  Larese}]{Chandra2021}
\bibinfo{author}{B.~Chandra}, \bibinfo{author}{V.~Singer},
  \bibinfo{author}{T.~Teschemacher}, \bibinfo{author}{R.~W\"{u}cheer},
  \bibinfo{author}{A.~Larese}, \bibinfo{title}{Nonconforming dirichlet boundary
  conditions in implicit material point method by means of penalty
  augmentation} \bibinfo{volume}{16} (\bibinfo{year}{2021})
  \bibinfo{pages}{2315--2335}.
\bibitem[{Charlton et~al.(2017)Charlton, Coombs and Augarde}]{Charlton2017}
\bibinfo{author}{T.J. Charlton}, \bibinfo{author}{W.M. Coombs},
  \bibinfo{author}{C.E. Augarde}, \bibinfo{title}{{iGIMP}: An implicit
  generalised interpolation material point method for large deformations},
  \bibinfo{journal}{Computers \& Structures} \bibinfo{volume}{190}
  (\bibinfo{year}{2017}) \bibinfo{pages}{108--125}.
\bibitem[{Coombs(2011)}]{coombs2011finitethesis}
\bibinfo{author}{W.~Coombs}, \bibinfo{title}{Finite deformation of particulate
  geomaterials: frictional and anisotropic Critical State elasto-plasticity},
  Ph.D. thesis, Durham University, \bibinfo{year}{2011}.
\bibitem[{Coombs and Augarde(2020)}]{Coombs2020}
\bibinfo{author}{W.M. Coombs}, \bibinfo{author}{C.E. Augarde},
  \bibinfo{title}{{AMPLE}: {A} {Material} {Point} {Learning} {Environment}},
  \bibinfo{journal}{Advances in Engineering Software} \bibinfo{volume}{139}
  (\bibinfo{year}{2020}) \bibinfo{pages}{102748}.
\bibitem[{Coombs et~al.(2020)Coombs, Augarde, Brennan, Brown, Charlton,
  Knappett, {Ghaffari Motlagh} and Wang}]{Coombs2020a}
\bibinfo{author}{W.M. Coombs}, \bibinfo{author}{C.E. Augarde},
  \bibinfo{author}{A.J. Brennan}, \bibinfo{author}{M.J. Brown},
  \bibinfo{author}{T.J. Charlton}, \bibinfo{author}{J.A. Knappett},
  \bibinfo{author}{Y.~{Ghaffari Motlagh}}, \bibinfo{author}{L.~Wang},
  \bibinfo{title}{On {L}agrangian mechanics and the implicit material point
  method for large deformation elasto-plasticity}, \bibinfo{journal}{Computer
  Methods in Applied Mechanics and Engineering} \bibinfo{volume}{358}
  (\bibinfo{year}{2020}) \bibinfo{pages}{112622}.
\bibitem[{Coombs et~al.(2018)Coombs, Charlton, Cortis and Augarde}]{Coombs2018}
\bibinfo{author}{W.M. Coombs}, \bibinfo{author}{T.J. Charlton},
  \bibinfo{author}{M.~Cortis}, \bibinfo{author}{C.E. Augarde},
  \bibinfo{title}{Overcoming volumetric locking in material point methods},
  \bibinfo{journal}{Compter Methods in Applied Mechanics and Engineering}
  \bibinfo{volume}{333} (\bibinfo{year}{2018}) \bibinfo{pages}{1--21}.
\bibitem[{Cortis et~al.(2018)Cortis, Coombs, Augarde, Brown, Brennan and
  Robinson}]{Cortis2018}
\bibinfo{author}{M.~Cortis}, \bibinfo{author}{W.M. Coombs},
  \bibinfo{author}{C.E. Augarde}, \bibinfo{author}{M.J. Brown},
  \bibinfo{author}{A.~Brennan}, \bibinfo{author}{S.~Robinson},
  \bibinfo{title}{Imposition of essential boundary conditions in the material
  point method}, \bibinfo{journal}{International Journal for Numerical Methods
  in Engineering} \bibinfo{volume}{113} (\bibinfo{year}{2018})
  \bibinfo{pages}{130--152}.
\bibitem[{{González Acosta} et~al.(2021){González Acosta}, Vardon and
  Hicks}]{Acosta2021}
\bibinfo{author}{J.L. {González Acosta}}, \bibinfo{author}{P.J. Vardon},
  \bibinfo{author}{M.A. Hicks}, \bibinfo{title}{Development of an implicit
  contact technique for the material point method}, \bibinfo{journal}{Computers
  and Geotechnics} \bibinfo{volume}{130} (\bibinfo{year}{2021})
  \bibinfo{pages}{103859}.
\bibitem[{Hammerquist and Nairn(2017)}]{Hammerquist2017}
\bibinfo{author}{C.C. Hammerquist}, \bibinfo{author}{J.A. Nairn},
  \bibinfo{title}{A new method for material point method particle updates that
  reduces noise and enhances stability}, \bibinfo{journal}{Computer Methods in
  Applied Mechanics and Engineering} \bibinfo{volume}{318}
  (\bibinfo{year}{2017}) \bibinfo{pages}{724--738}.
\bibitem[{Hansbo et~al.(2017)Hansbo, Larson and Larsson}]{Hansbo2017}
\bibinfo{author}{P.~Hansbo}, \bibinfo{author}{M.~Larson},
  \bibinfo{author}{K.~Larsson}, \bibinfo{title}{Cut finite element methods for
  linear elasticity problems}, \bibinfo{journal}{Lecture Notes in Computational
  Science and Engineering}  (\bibinfo{year}{2017}) \bibinfo{pages}{25--63}.
\bibitem[{Harlow(1964)}]{Harlow1964}
\bibinfo{author}{F.~Harlow}, \bibinfo{title}{The particle-in-cell computing
  method for fluid dynamics}, \bibinfo{journal}{Methods for Computational
  Physics} \bibinfo{volume}{3} (\bibinfo{year}{1964}).
\bibitem[{Love and Sulsky(2006)}]{Love2006}
\bibinfo{author}{E.~Love}, \bibinfo{author}{D.L. Sulsky}, \bibinfo{title}{An
  energy-consistent material-point method for dynamic finite deformation
  plasticity}, \bibinfo{journal}{International Journal for Numerical Methods in
  Engineering} \bibinfo{volume}{65} (\bibinfo{year}{2006})
  \bibinfo{pages}{1608--1638}.
\bibitem[{Ma et~al.(2010)Ma, Giguere, Jayaraman and Zhang}]{Ma2010}
\bibinfo{author}{X.~Ma}, \bibinfo{author}{P.T. Giguere},
  \bibinfo{author}{B.~Jayaraman}, \bibinfo{author}{D.Z. Zhang},
  \bibinfo{title}{Distribution coefficient algorithm for small mass nodes in
  material point method}, \bibinfo{journal}{Journal of Computational Physics}
  \bibinfo{volume}{229} (\bibinfo{year}{2010}) \bibinfo{pages}{7819 -- 7833}.
\bibitem[{Molstad(1977)}]{Molstad1977}
\bibinfo{author}{T.~Molstad}, \bibinfo{title}{Finite deformation analysis using
  the finite element method}, Ph.D. thesis, University of British Columbia,
  \bibinfo{year}{1977}.
\bibitem[{Nairn and Hammerquist(2021)}]{Nairn2021}
\bibinfo{author}{J.A. Nairn}, \bibinfo{author}{C.C. Hammerquist},
  \bibinfo{title}{Material point method simulations using an approximate full
  mass matrix inverse}, \bibinfo{journal}{Computer Methods in Applied Mechanics
  and Engineering} \bibinfo{volume}{377} (\bibinfo{year}{2021})
  \bibinfo{pages}{113667}.
\bibitem[{Pretti et~al.(2022)Pretti, Coombs, Augarde, Sims, Puigvert and
  Guti\'{e}rrez}]{Pretti2022}
\bibinfo{author}{G.~Pretti}, \bibinfo{author}{W.M. Coombs},
  \bibinfo{author}{C.E. Augarde}, \bibinfo{author}{B.~Sims},
  \bibinfo{author}{M.M. Puigvert}, \bibinfo{author}{J.A.R. Guti\'{e}rrez},
  \bibinfo{title}{An updated lagrangian, energy conserving material point
  method for dynamic analysis}, \bibinfo{journal}{Computer Methods in Applied
  Mechanics and Engineering (under review)}  (\bibinfo{year}{2022}).
\bibitem[{Remmerswaal(2017)}]{Remmerswaal2017}
\bibinfo{author}{G.~Remmerswaal}, \bibinfo{title}{Development and
  implementation of moving boundary conditions in the Material Point Method},
  Master's thesis, TU Delft, \bibinfo{year}{2017}.
\bibitem[{Sadeghirad et~al.(2013)Sadeghirad, Brannon and
  Guilkey}]{Sadeghirad2013a}
\bibinfo{author}{A.~Sadeghirad}, \bibinfo{author}{R.~Brannon},
  \bibinfo{author}{J.~Guilkey}, \bibinfo{title}{Second-order convected particle
  domain interpolation {(CPDI2)} with enrichment for weak discontinuities at
  material interfaces}, \bibinfo{journal}{International Journal for Numerical
  Methods in Engineering} \bibinfo{volume}{95} (\bibinfo{year}{2013})
  \bibinfo{pages}{928--952}.
\bibitem[{Sadeghirad et~al.(2011)Sadeghirad, Brannon and
  Burghardt}]{Sadeghirad2011}
\bibinfo{author}{A.~Sadeghirad}, \bibinfo{author}{R.M. Brannon},
  \bibinfo{author}{J.~Burghardt}, \bibinfo{title}{A convected particle domain
  interpolation technique to extend applicability of the material point method
  for problems involving massive deformations}, \bibinfo{journal}{International
  Journal for Numerical Methods in Engineering} \bibinfo{volume}{86}
  (\bibinfo{year}{2011}) \bibinfo{pages}{1435--1456}.
\bibitem[{Simo(1992)}]{Simo1992a}
\bibinfo{author}{J.~Simo}, \bibinfo{title}{Algorithms for static and dynamic
  multiplicative plasticity that preserve the classical return mapping schemes
  of the infinitesimal theory}, \bibinfo{journal}{Computer Methods in Applied
  Mechanics and Engineering} \bibinfo{volume}{99} (\bibinfo{year}{1992})
  \bibinfo{pages}{61--112}.
\bibitem[{Solowski et~al.(2021)Solowski, Berzins, Coombs, Guilkey, Möller,
  Tran, Adibaskoro, Seyedan, Tielen and Soga}]{Solowski2021}
\bibinfo{author}{W.T. Solowski}, \bibinfo{author}{M.~Berzins},
  \bibinfo{author}{W.M. Coombs}, \bibinfo{author}{J.E. Guilkey},
  \bibinfo{author}{M.~Möller}, \bibinfo{author}{Q.A. Tran},
  \bibinfo{author}{T.~Adibaskoro}, \bibinfo{author}{S.~Seyedan},
  \bibinfo{author}{R.~Tielen}, \bibinfo{author}{K.~Soga},
  \bibinfo{title}{Material point method: Overview and challenges ahead},
  volume~\bibinfo{volume}{54} of \textit{\bibinfo{series}{Advances in Applied
  Mechanics}}, \bibinfo{publisher}{Elsevier}, \bibinfo{year}{2021}, pp.
  \bibinfo{pages}{113--204}.
\bibitem[{de~Souza~Neto et~al.(2008)de~Souza~Neto, Peric and Owen}]{SouzaNeto}
\bibinfo{author}{E.A. de~Souza~Neto}, \bibinfo{author}{D.~Peric},
  \bibinfo{author}{D.R.J. Owen}, \bibinfo{title}{Computational Methods For
  Plasticity: Theory and Applications}, \bibinfo{publisher}{John Wiley \& Sons,
  Ltd}, \bibinfo{year}{2008}.
\bibitem[{Steffen et~al.(2008)Steffen, Wallstedt, Guilkey, Kirby and
  Berzins}]{Steffen2008}
\bibinfo{author}{M.~Steffen}, \bibinfo{author}{P.C. Wallstedt},
  \bibinfo{author}{J.E. Guilkey}, \bibinfo{author}{R.M. Kirby},
  \bibinfo{author}{M.~Berzins}, \bibinfo{title}{Examination and analysis of
  implementation choices within the material point method},
  \bibinfo{journal}{Computer Modeling in Engineering and Sciences}
  \bibinfo{volume}{31} (\bibinfo{year}{2008}) \bibinfo{pages}{107--127}.
\bibitem[{Sticko et~al.(2020)Sticko, Ludvigsson and Kreiss}]{Sticko2020}
\bibinfo{author}{S.~Sticko}, \bibinfo{author}{G.~Ludvigsson},
  \bibinfo{author}{G.~Kreiss}, \bibinfo{title}{High-order cut finite elements
  for the elastic wave equation}, \bibinfo{journal}{Advances in Computational
  Mathematics} \bibinfo{volume}{46:45} (\bibinfo{year}{2020}) \bibinfo{pages}{1
  -- 28}.
\bibitem[{Sulsky et~al.(1994)Sulsky, Chen and Schreyer}]{Sulsky1994}
\bibinfo{author}{D.~Sulsky}, \bibinfo{author}{Z.~Chen}, \bibinfo{author}{H.L.
  Schreyer}, \bibinfo{title}{A particle method for history-dependent
  materials}, \bibinfo{journal}{Computer Methods in Applied Mechanics and
  Engineering} \bibinfo{volume}{118} (\bibinfo{year}{1994})
  \bibinfo{pages}{179--196}.
\bibitem[{de~Vaucorbeil et~al.(2020)de~Vaucorbeil, Nguyen, Sinaie and
  Wu}]{Vaucorbeil2020}
\bibinfo{author}{A.~de~Vaucorbeil}, \bibinfo{author}{V.P. Nguyen},
  \bibinfo{author}{S.~Sinaie}, \bibinfo{author}{J.Y. Wu},
  \bibinfo{title}{Material point method after 25 years: Theory, implementation,
  and applications}, \bibinfo{journal}{Advances in applied mechanics}
  \bibinfo{volume}{53} (\bibinfo{year}{2020}) \bibinfo{pages}{185--398}.
\bibitem[{Wang et~al.(2016)Wang, Vardon and Hicks}]{Wang2016}
\bibinfo{author}{B.~Wang}, \bibinfo{author}{P.~Vardon},
  \bibinfo{author}{M.~Hicks}, \bibinfo{title}{Investigation of retrogressive
  and progressive slope failure mechanisms using the material point method},
  \bibinfo{journal}{Computers and Geotechnics} \bibinfo{volume}{78}
  (\bibinfo{year}{2016}) \bibinfo{pages}{88 -- 98}.
\bibitem[{Wang et~al.(2019)Wang, Coombs, Augarde, Cortis, Charlton, Brown,
  Knappett, Brennan, Davidson, Richards and Blake}]{Wang2019}
\bibinfo{author}{L.~Wang}, \bibinfo{author}{W.~Coombs},
  \bibinfo{author}{C.~Augarde}, \bibinfo{author}{M.~Cortis},
  \bibinfo{author}{T.~Charlton}, \bibinfo{author}{M.~Brown},
  \bibinfo{author}{J.~Knappett}, \bibinfo{author}{A.~Brennan},
  \bibinfo{author}{C.~Davidson}, \bibinfo{author}{D.~Richards},
  \bibinfo{author}{A.~Blake}, \bibinfo{title}{On the use of domain-based
  material point methods for problems involving large distortion},
  \bibinfo{journal}{Computer Methods in Applied Mechanics and Engineering}
  \bibinfo{volume}{355} (\bibinfo{year}{2019}) \bibinfo{pages}{1003--1025}.
\bibitem[{Wang et~al.(2021)Wang, Coombs, Augarde, Cortis, Brown, Brennan,
  Knappett, Davidson, Richards, White and Blake}]{Wang2021}
\bibinfo{author}{L.~Wang}, \bibinfo{author}{W.M. Coombs}, \bibinfo{author}{C.E.
  Augarde}, \bibinfo{author}{M.~Cortis}, \bibinfo{author}{M.J. Brown},
  \bibinfo{author}{A.J. Brennan}, \bibinfo{author}{J.A. Knappett},
  \bibinfo{author}{C.~Davidson}, \bibinfo{author}{D.~Richards},
  \bibinfo{author}{D.J. White}, \bibinfo{author}{A.P. Blake},
  \bibinfo{title}{An efficient and locking-free material point method for
  three-dimensional analysis with simplex elements},
  \bibinfo{journal}{International Journal for Numerical Methods in Engineering}
  \bibinfo{volume}{122} (\bibinfo{year}{2021}) \bibinfo{pages}{3876--3899}.

\end{thebibliography}

%

\end{document}